\documentclass[review,sort&compress,1p,times]{elsarticle}

\usepackage{makecell}
\usepackage{amsmath}
\usepackage{pdfpages}
\usepackage{geometry}  
\usepackage[colorlinks,linkcolor=green]{hyperref}   
\usepackage{float}
\usepackage{amsfonts,amssymb}
\usepackage{bm}  
\usepackage{graphicx}
\usepackage{mathrsfs}
\usepackage{cancel} 
\usepackage{extarrows} 
\usepackage{subeqnarray}
\usepackage{subcaption}
\usepackage{cases}   
\usepackage{xcolor} 
\usepackage{colortbl,booktabs} 
\usepackage{multirow}
\usepackage{threeparttable}  
\usepackage{caption}  
\usepackage{hyperref}
\usepackage{algorithm}
\newcommand{\ud}{\,\mathrm{d}}

\journal{Engineering with Computers}

\begin{document}

\begin{frontmatter}
\title{Virtual element method for thermomechanical analysis of electronic packaging structures with multi-scale features}

\author[address1]{Yanpeng Gong}
\ead{yanpeng.gong@bjut.edu.cn}

\author[address1]{Sishuai Li}

\author[address1]{Fei Qin}

\author[address2]{Bingbing Xu\corref{mycorrespondingauthor}}
\ead{bingbing.xu@ikm.uni-hannover.de}

\address[address1]{Institute of Electronics Packaging Technology and Reliability, Department of Mechanics, Beijing University of Technology, Beijing, 100124, China}
\address[address2]{Institute of Continuum Mechanics, Leibniz University Hannover, Hannover, Germany}
\cortext[mycorrespondingauthor]{Corresponding author}

\begin{abstract}
This paper presents two approaches: the virtual element method (VEM) and the stabilization-free virtual element method (SFVEM) for analyzing thermomechanical behavior in electronic packaging structures with geometric multi-scale features.
Since the virtual element method allows the use of arbitrary polygonal elements,
the inherent mesh flexibility of VEM allows localized mesh modifications without affecting global mesh structure, 
making it particularly effective for the analysis of electronic packaging reliability involving complex geometries and multiple geometric scales.
The approach implements a novel non-matching mesh generation strategy that strategically combines polygonal meshes for complex small-scale regions with regular quadrilateral meshes for larger domains. 
The VEM formulation addresses both heat conduction and thermomechanical coupling problems, with comprehensive verification through analytical benchmarks and practical electronic packaging case studies, including Through-Silicon Via (TSV), Ball Grid Array (BGA), and Plastic Ball Grid Array (PBGA) structures.
Results demonstrate that the method accurately captures stress concentrations at material interfaces and provides reliable thermal and mechanical response predictions. 
Some MATLAB codes for the numerical examples are provided at https://github.com/yanpeng-gong/VEM-electronic-packaging and on the VEMhub website (www.vemhub.com).
\end{abstract}


\end{frontmatter}

\section{Introduction}
\label{sec1}

Electronic packaging transforms semiconductor devices into functional products, critically impacting performance and reliability throughout production and assembly~\cite{moreau2012semiconductor}.
Since Moore's Law emerged in 1965, it has driven semiconductor advancement, but as technology approaches physical limits with sub-3 nm mobile phone chips, 
packaging must evolve toward miniaturization and higher density integration~\cite{ahmed2020interfacial}. 
Modern 3D integrated circuits face complex thermal and mechanical loads~\cite{RN1745,gong2022isogeometric}, 
which cause over $55\%$ of electronic device failures~\cite{he2021thermal}. 
Consequently, heat conduction and mechanical analysis have become essential for ensuring electronic packaging reliability.

Numerical analysis offers advantages of simplicity, cost-effectiveness, and accuracy in electronic packaging reliability assessments. 
Researchers have developed diverse simulation methods including Finite Element Method (FEM)~\cite{VENKATANAGACHANDANA2022675}, 
Boundary Element Method~\cite{yu2021thermal,gong2022isogeometric}, FEM/BEM coupling~\cite{RN1793,RN1014}, and Finite Volume Method (FVM)~\cite{khor2010fvm} to address thermal and mechanical challenges. 
Oukaira et al.~\cite{oukaira2023thermal} combined Field-Programmable Gate Array with FEM for transient thermal analysis of System-in-Package models. 
Gong et al.~\cite{gong2022isogeometric} introduced Isogeometric Boundary Element Method for analyzing heat transfer in multi-scale structures with arbitrary heat sources. 
Feng et al.~\cite{feng2023engineered} developed Stable Node-based Smoothed FEM for electro-thermal-mechanical coupling in IGBT modules, 
while other researchers proposed thermal-mechanical coupled phase field models for interconnect structure analysis~\cite{RN1745}. 
Yao et al.~\cite{yao2022physics} applied physics-based nested artificial neural networks to assess warpage reliability in fan-out wafer-level packaging, 
and FEM-BEM coupling methods have been employed to examine various electronic packaging problems~\cite{RN1793,RN1014}. Compared to FEM and other conventional methods, the Virtual Element Method can utilize arbitrary polygonal elements, making it particularly suitable for analyzing models with complex geometric structures. This method exhibits excellent geometric adaptability with lower mesh quality requirements, maintaining algorithmic convergence even under mesh distortion conditions. Therefore, VEM demonstrates superior geometric robustness and excellent convergence properties for any element shape, including non-convex elements. For geometric multi-scale structures, VEM enables local mesh refinement in specialized regions using non-matching meshes without affecting the overall mesh distribution, thereby reducing element count and computational cost. When analyzing models with extensive geometric multi-scale features using FEM, numerous transition elements are typically required, leading to high computational costs. In contrast, VEM can employ arbitrary-shaped elements to replace transition regions, significantly reducing the number of transition elements. Consequently, this method demonstrates tremendous potential for addressing such problems.

Despite widespread application of existing methods in electronic packaging reliability analysis, 
geometric multi-scale structures present persistent challenges. 
FEM simulations require rational element division to avoid distortion, but significant size differences necessitate transition meshes-often resulting in millions of elements that increase computational burden and compromise convergence. 
Additionally, many methods have limited applications; BEM excels with infinite domains or homogeneous problems, 
While FVM is applicable to both solid mechanics and fluid flow analysis, it is more widely used in fluid mechanics.
The Virtual Element Method (VEM) is an innovative FEM extension which offers compelling advantages for these challenges \cite{liguori2024hybrid, xu2024high, beirao2013basic}. 
VEM's mesh flexibility accommodates polygons of arbitrary shape, simplifying complex geometry discretization. 
Its non-matching mesh capability enables targeted local refinement without disrupting overall mesh distribution, reducing element count and computational demands. 
Furthermore, by eliminating the need for explicit basis function expressions, VEM enhances both accuracy and efficiency-making it particularly valuable for reliability analysis of geometric multi-scale electronic packaging models.

VEM was first introduced by Beirao da Veiga et al.~\cite{beirao2013basic} for linear elasticity problems, and has since gained traction in various fields.
Artioli et al.~\cite{artioli2017arbitrary} demonstrated VEM's accuracy in two-dimensional linear elasticity by comparing displacement results with FEM on structures. 
Subsequent studies confirmed VEM's geometric robustness and excellent convergence properties regardless of element shape~\cite{sorgente2022polyhedral,van2021virtual}. 
VEM applications have expanded beyond linear elasticity to nonlinear mechanical problems.
Liu et al.~\cite{liu2023virtual} developed VEM for phase-field modeling of dynamic fracture using explicit time integration, 
while Hussein et al.~\cite{hussein2019computational} created an efficient low-order VEM for crack propagation by splitting elements along growth directions. 
Cihan et al.~\cite{cihan2022virtual} applied VEM to complex contact problems using node-to-node discretization to avoid mesh matching complications, 
and others implemented VEM for thermomechanical coupling analysis through Abaqus UEL tools~\cite{dhanush2019implementation}. Additionally, Aldakheel et al.~\cite{Aldakheel2019} developed a low-order virtual element scheme for finite strain thermo-plasticity problems, demonstrating the method's capability in handling complex nonlinear thermomechanical coupling scenarios.

The basic idea of VEM is to construct a projection $\Pi_k^\nabla$ for the displacement and temperature fields, and a remainder term $\bm{u}-\Pi_k^\nabla\bm{u}$ is needed.
The remainder leads to additional stabilization terms in the VEM formulation.
Of course, for multi-field coupling problems, we often do not want any stabilization terms.
This is because stabilization terms contain user-defined parameters whose selection affects the computation of both temperature and displacement fields. In multi-field coupling problems, these parameter choices may cause error propagation and amplification between different physical fields, while different stabilization formulations may lead to variations in results and increased computational complexity.
In this work, we propose a VEM formulation for thermomechanical coupling problems without any stabilization.
The basic idea of the stabilization-free VEM is to construct a higher-order projection for the gradient.
By choosing the appropriate polynomial order, we can obtain a stiffness matrix that does not require any stabilization term.
The stabilization-free VEM has been successfully applied to various problems, including linear elasticity \cite{Xu2024} and nonlinear elasticity \cite{xu2024stabilization,Xu2023,Xu2024a}.
But up to now, the stabilization-free VEM has not been applied to thermomechanical coupling problems.

Besides, this work leverages VEM's geometric mesh flexibility and automatic generation capabilities to address reliability issues in electronic packaging structures with substantial scale differences. 
By applying polygonal meshes to small-scale regions requiring refined analysis while using coarser meshes for larger areas, and integrating these through non-matching mesh techniques, 
we achieve high-precision local simulation while minimizing element count and computational demands. 
In this work, we investigate both linear elastic and thermomechanical coupling problems based on conventional VEM and stabilization-free VEM to analyze mechanical performance of these structures. 
We demonstrate VEM's effectiveness and accuracy through numerical examples, including analytical models and practical electronic packaging structures. 
Our results confirm VEM's potential for efficient and accurate thermal and mechanical analysis in electronic packaging structures with multi-scale features.
Relevant codes for VEM in thermomechanical coupling problems are available on GitHub (https://github.com/yanpeng-gong/VEM-electronic-packaging) and the VEMhub website (www.vemhub.com).

The paper is structured as follows: Section \ref{s2} introduces temperature and force field governing equations using VEM, 
while Section \ref{s3} details VEM formulation for thermoelastic problems including heat conduction and thermomechanical coupling. 
The idea of the stabilization-free VEM and the formulation for the thermomechanical coupling is given in Section \ref{s4}.
Section \ref{sec:meshing} describes our non-matching mesh generation strategy. 
Section \ref{Examples} applies VEM to numerical examples analyzing mechanical performance in analytical models and electronic packaging structures. 
Section \ref{conclusions} concludes with results and discussion.

\section{Governing equations for VEM in thermoelastic problems}
\label{s2}
\begin{figure}[htbp]
    \centering
    \includegraphics[width=0.35\textwidth]{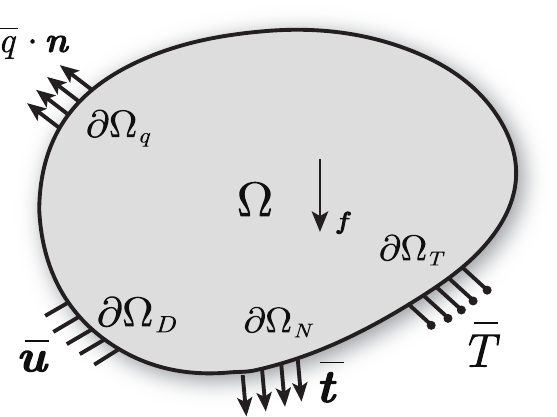}
 \caption{   Schematic diagram of the thermoelastic problem domain.}
    \label{fig:thermoelastic problem domain}
\end{figure}

The thermoelastic problem is defined on continuous domain $\Omega$ (Fig.~\ref{fig:thermoelastic problem domain}) 
with body force $f$ per unit area. 
In thermomechanical coupling analysis, the domain boundary $\partial \Omega$ is partitioned differently for each physical field: the thermal conduction problem divides the boundary into two disjoint regions $\partial \Omega=\partial \Omega_T \cup \partial \Omega_q$, while the elastic problem similarly divides the boundary into two disjoint regions $\partial \Omega=\partial \Omega_D \cup \partial \Omega_N$. As thermomechanical coupling involves both temperature and force fields, governing equations for both physical fields must be established.

\subsection{Steady-state heat conduction governing equations}

For the steady-state heat conduction problem, the governing equation in domain $\Omega$ is
\begin{equation}\label{eq:heat_conduction}
    \nabla\left(\lambda_k\nabla T\right) +Q = 0,
\end{equation}
where $\lambda_k$ is the thermal conductivity and $T$ is the temperature field, $Q$ is the internal heat generation per unit volume.

The boundary conditions are prescribed on $\partial\Omega = \partial\Omega_T \cup \partial\Omega_q$ as
\begin{equation}
    T(\bm{x}) = \overline{T}, \quad \forall \bm{x}\text{ on }\partial\Omega_T,
\end{equation}
\begin{equation}
    -\lambda_k\nabla T(\bm{x})\cdot\bm{n} = \overline{q}, \quad \forall \bm{x}\text{ on }\partial\Omega_q,
\end{equation}
where $\overline{q}$ denotes the prescribed heat flux, and $\bm{n}$ is the outward unit normal vector.


The weak formulation of the thermal equilibrium equation is 
\begin{equation}
    \label{s2.bilinearT}
      \lambda_k a_T(T,\varphi)+ \int_\Omega Q\varphi\mathrm{~d}\Omega= \ell_T(\varphi),
\end{equation}
where $\varphi$ is variational of temperature. For simplicity, we assume $Q = 0$ (no internal heat generation) in this study.

The associated bilinear and linear forms are defined as
\begin{equation}
    a_T(T,\varphi) = \int_\Omega \nabla T \cdot \nabla \varphi \,\ud\Omega,
\end{equation}
\begin{equation}
    \ell_T(\varphi) = \int_{\partial\Omega_q} \varphi \overline{q} \,\ud\Gamma.
\end{equation}

\subsection{Thermomechanical coupling governing equations}
\label{sec:thermomechanical governing equations}
For the thermoelastic problem, the governing equation of motion in domain $\Omega$ takes the form
\begin{equation}\label{eq:equilibrium_equation}
    \nabla \cdot \bm{\sigma} + \bm{f} = \bm{0},
\end{equation}
where $\bm{\sigma}$ is the Cauchy stress tensor and $\bm{f}$ denotes the body force vector. 
The constitutive relationship considering thermal strain is given by
\begin{equation}\label{eq:constitutive_equation}
    \bm{\sigma} = \mathbb{D} : (\bm{\varepsilon} - \bm{\varepsilon}_t),
\end{equation}
with the strain tensor is being defined as
\begin{equation}\label{eq:strain_displacement}
    \bm{\varepsilon}(\bm{u}) = \frac{1}{2}(\nabla\bm{u} + \nabla\bm{u}^T).
\end{equation}

The thermal strain is expressed as
\begin{equation}\label{eq:thermal_strain}
    \boldsymbol{\varepsilon}_t=\alpha\left(T-T_{ref}\right)\boldsymbol{I},
\end{equation}
where $T_{ref}$ represents the reference (initial) temperature at which thermal strains are zero. This formulation correctly captures that thermal strain results from temperature changes relative to a reference state, not from absolute temperature values. For instance, maintaining a constant temperature of 298 K would indeed produce no thermal strain if $T_{ref} = 298$ K, which is physically consistent. In addition, $\bm{I}$ is the identity tensor,
$\alpha$ is the thermal expansion coefficient, and $\mathbb{D}$ is the fourth-order elasticity tensor, which can be obtained from the energy density function $\Psi$ as
\begin{equation}
    \Psi(\bm{\varepsilon}) = \frac{\lambda}{2}(\text{tr}\bm{\varepsilon})^2+\mu\bm{\varepsilon}:\bm{\varepsilon},
\end{equation}
where $\lambda$ and $\mu$ are the $\mathrm{Lam\acute{e}}$ parameters.

The elasticity tensor $\mathbb{D}$ is derived from the energy density function as \cite{BERBATOV2021351}
\begin{equation}
D_{ijkl}=\frac{\partial^{2}\Psi}{\partial\varepsilon_{ij}\partial\varepsilon_{kl}}=\lambda\delta_{ij}\delta_{kl}+\mu\left(\delta_{ik}\delta_{jl}+\delta_{il}\delta_{jk}\right).
\end{equation}
The boundary conditions are prescribed on $\partial\Omega = \partial\Omega_D \cup \partial\Omega_N$ as
\begin{align}
    \bm{u} &= \overline{\bm{u}}, \quad \text{on }\partial\Omega_D \label{eq:dirichlet_bc} \\[8pt]
    \bm{\sigma} \cdot \bm{n} &= \overline{\bm{t}}, \quad \text{on }\partial\Omega_N \label{eq:neumann_bc}
\end{align}
where $\overline{\bm{u}}$ and $\overline{\bm{t}}$ represent the prescribed displacement and traction vector, respectively, and $\bm{n}$ denotes the outward unit normal vector.

The weak formulation of the equilibrium equation reads: find $\bm{u} \in \gamma_u$ such that
\begin{equation}\label{eq:weak_form_mechanical}
    a_u(\bm{u},\bm{v}) = \ell_u(\bm{v}), \quad \forall \bm{v} \in \gamma_u
\end{equation}
where $\bm{v}$ is the test function and $\gamma_u = \{\bm{u} \in [\mathcal{H}^1(\Omega)]^2: \bm{u}|_{\partial\Omega_D} = \overline{\bm{u}}\}$ is the displacement solution space. Both the trial function $\bm{u}$ and test function $\bm{v}$ belong to the same displacement solution space $\gamma_u$.

The associated bilinear and linear forms are defined as \cite{mengolini2019engineering}
\begin{equation}
    \label{eq:bilinear_form_mechanical}
    a_u(\bm{u},\bm{v}) = \int_\Omega \bm{\varepsilon}(\bm{u}): \mathbb{D} : \bm{\varepsilon}(\bm{v}) \,\ud\Omega ,
\end{equation}
\begin{equation}
    \label{eq:linear_form_mechanical}
    \ell_u(\bm{v}) = \int_\Omega \bm{v} \cdot \bm{f} \,d\Omega + \int_{\partial\Omega_N} \bm{v} \cdot \overline{\bm{t}} \,\ud\Gamma 
    +\int_\Omega \bm{\varepsilon}(\bm{v}):\mathbb{D}:\bm{\varepsilon}_t\ud\Omega.
\end{equation}
The last term in Eq. \eqref{eq:linear_form_mechanical} is the body force generated by thermal expansion.

Generally speaking, the solution for steady-state thermomechanical coupling problems is to first solve the temperature field 
and then apply the temperature as a load to the structural analysis. For the electronic packaging structures considered in this work, material properties do not exhibit significant temperature dependence within the operating range, allowing us to employ this decoupled approach without loss of accuracy. By solving the thermal and mechanical problems sequentially, we can utilize specialized solvers optimized for each individual field, which enhances numerical stability and computational efficiency compared to solving the fully coupled system. 
This approach reduces the computational cost compared to solving the fully coupled system.
Therefore, we need to construct discrete formats for the bilinear formats of the temperature field and displacement field respectively.

\section{VEM formulation for thermomechanical coupling}
\label{s3}
\begin{figure}[htbp]
    \centering
    \begin{subfigure}[b]{0.7\textwidth}
        \centering
    \includegraphics[width=0.8\textwidth]{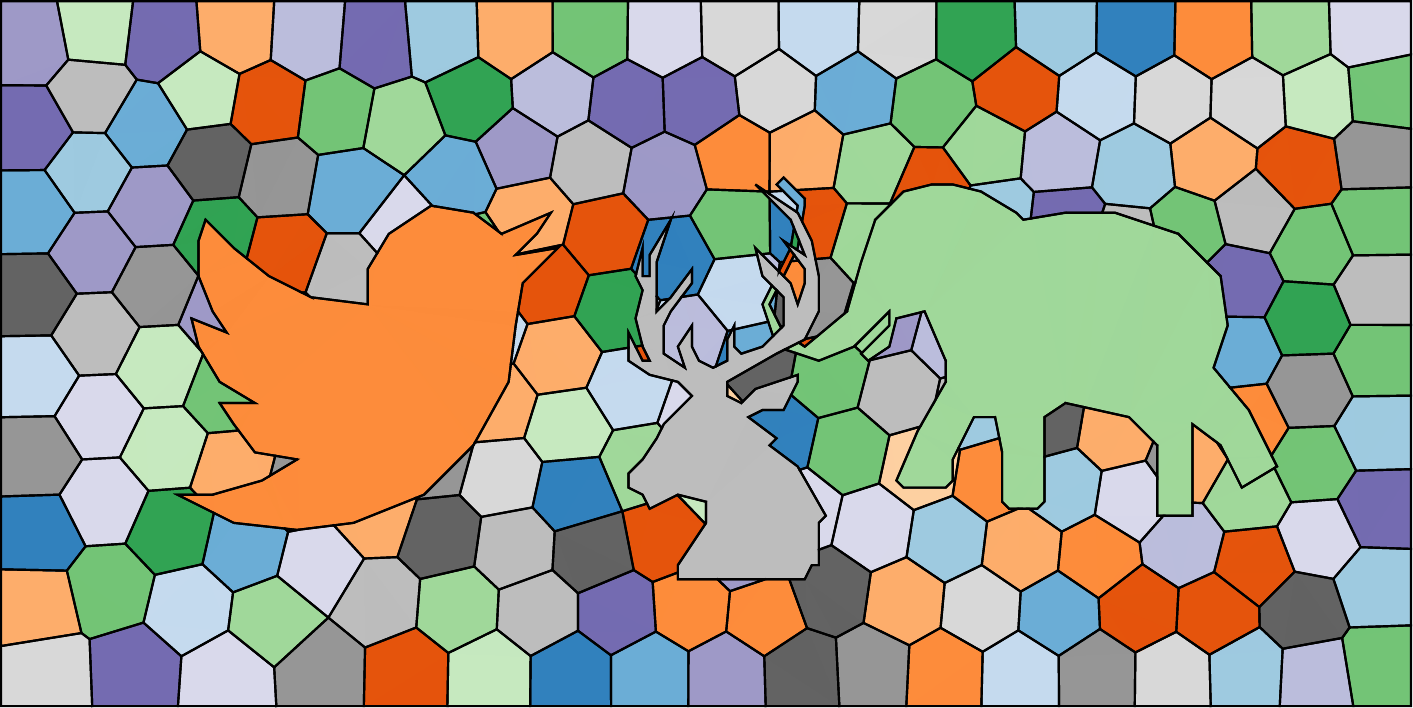}
    \caption{ Polygonal mesh discretization with arbitrary element shapes including irregular polygons, animal-like shapes (bird, deer and elephant), pentagons and hexagons.}
    \label{s3.mesh}
    \end{subfigure}
    \hfill
    
    \begin{subfigure}[b]{0.75\textwidth}
        \centering
    \includegraphics[width=0.85\textwidth]{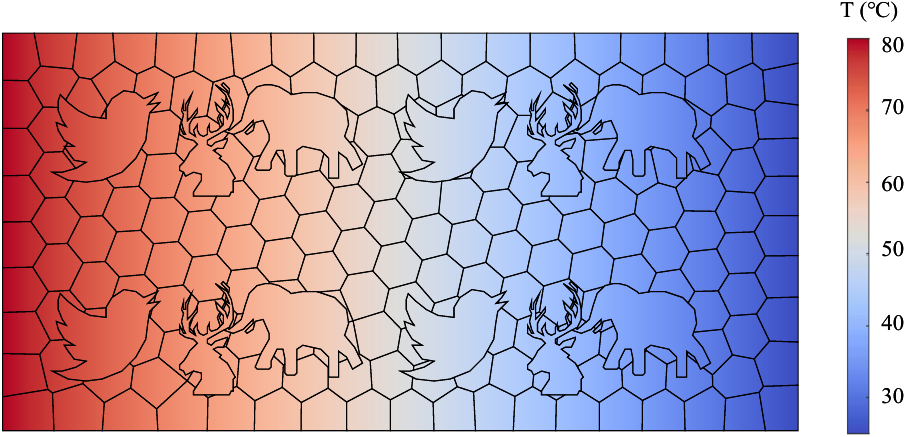}
   \caption{The temperature field of the computational domain.}
    \label{s4.temperature field}
    \end{subfigure}
 \caption{Polygonal mesh and the temperature field for the computational domain.}
\end{figure}

The VEM requires the discretization of the computational domain into finite subdomains, 
such that $\Omega^h \equiv \text{int}\left(\bigcup_{i} E_i\right)$, where the operator `int' denotes the interior of a set, and  $\Omega^h $ represents the union of the interior regions of all elements $ E_i$, excluding the boundaries,
where each element $E_i$ can be a polygon with arbitrary number of nodes and shape, as shown in Fig.\ref{s3.mesh}. Fig.\ref{s4.temperature field} illustrates both the irregular polygonal mesh capability and the computed temperature field from a representative heat conduction analysis, with boundary conditions of $80^{\circ} \mathrm{C}$ on the left boundary, $25^{\circ} \mathrm{C}$ on the right boundary, and adiabatic conditions on the top and bottom boundaries. This example demonstrates the mesh flexibility of VEM and its ability to produce smooth, physically reasonable temperature distributions on highly irregular mesh topologies.
In this work, the analysis process involves first solving the heat transfer problem to obtain the temperature field,
followed by solving the mechanical problem.

\subsection{VEM formulation for heat conduction} 
\label{sec:VEM_heat_conduction}

For the heat conduction problem, we first introduce the scalar virtual element function space
\begin{equation}\label{eq:VEM_space_thermal}
    \mathcal{V}_T(E) = \left\{T \in \mathcal{H}^1(E): T|_{\partial E} \in \mathcal{P}_k(\partial E), T = \overline{T} \text{ on }\partial\Omega_T\right\}.
\end{equation}
where $\mathcal{P}_k(E)$ denotes the space of polynomial functions of degree not exceeding $k$. 
In this work, we choose $k = 1$. 
In the polynomial space, the scaled coordinates $\xi$ and $\eta$ are defined as~\cite{artioli2017arbitrary}
\begin{equation}
    \label{eq:scaled_coordinates}
    \xi = \frac{x - \bar{x}}{h_E} \quad\eta = \frac{y - \bar{y}}{h_E}.
\end{equation}
where $\bar{\bm{x}} = (\bar{x},\bar{y})$ denotes the centroid of element $E$, and $h_E$ represents the maximum distance between any two nodes in element $E$.

Since the element shape functions do not have explicit expressions, 
we can define a projection operator to compute the bilinear form: 
mapping from the approximate function space $\mathcal{V}^h_T(E)$ on element $E$ to the polynomial space $\mathcal{P}_k(E)$, 
that is, $\Pi^\nabla: \mathcal{V}^h_T(E) \to \mathcal{P}_k(E)$:
\begin{equation}\label{eq:projection_orthogonality}
      a_T^E(\varphi^h - \Pi^\nabla\varphi^h,p) = 0, \quad \forall p \in \mathcal{P}_k(E).
\end{equation}

Expanding Eq.~\eqref{eq:projection_orthogonality} yields
\begin{equation}\label{eq:projection_expanded-thermal}
    \int_E \nabla(\Pi^\nabla\varphi^h) \cdot \nabla p \ \mathrm{~d} \Omega = \int_E \nabla\varphi^h \cdot \nabla p \mathrm{~d} \Omega.
\end{equation}

According to Eq.~\eqref{eq:projection_expanded-thermal}, it can be obtained
\begin{equation}\label{eq:projection_expanded2-thermal}
   \int_E \Pi^{\nabla} \varphi^h \mathrm{~d} \Omega=\int_E \varphi^h \mathrm{~d} \Omega.
\end{equation}

The projection operator can be solved by~\cite{dhanush2019implementation}
\begin{equation}\label{eq:projection_expanded2-thermal-1}
    \tilde{\Pi}=G^{-1} B,
\end{equation}
where
\begin{equation}\label{eq:projection_expanded2-thermal-G}
   \boldsymbol{G}=a_T^E\left(\boldsymbol{p}_\alpha, \boldsymbol{p}_\alpha\right)=\int_E \nabla \boldsymbol{p}_\alpha \cdot\lambda_k \nabla \boldsymbol{p}_\alpha \mathrm{~d} \Omega,
\end{equation}
\begin{equation}\label{eq:projection_expanded2-thermal-B}
   \boldsymbol{B}=a_T^E\left(\boldsymbol{\psi}_i, p_\alpha\right)=-\int_E \boldsymbol{\psi}_i \cdot\left[\nabla \cdot\left(\lambda_k\nabla \boldsymbol{p}_\alpha\right)\right] \mathrm{~d} \Omega+\int_{\partial E} \boldsymbol{\psi}_i \cdot \left(\lambda_k\nabla\boldsymbol{p}_\alpha\right) \boldsymbol{n} \mathrm{~d} \Gamma,
\end{equation}
where $\psi$ represents the basis functions of $\mathcal{V}_T(E)$. Since $\boldsymbol{p}_\alpha \in P_k(E)$, we have $\nabla \cdot\left(\lambda_k \nabla \boldsymbol{p}_\alpha\right) \in P_{k-2}(E)$. The first term on the right-hand side of Eq.~\eqref{eq:projection_expanded2-thermal-B} relates to internal degrees of freedom. For our polynomial order $k = 1$, this first term is not considered in the computation.

For temperature field $T^h$, the bilinear form on each element can be written as:
\begin{align}
\label{eq:bilinear_decomposition_thermal}
a_T^E(T^h,\varphi^h) &=a_T^E(T^h - \Pi^\nabla T^h + \Pi^\nabla T^h,\varphi^h - \Pi^\nabla\varphi^h + \Pi^\nabla\varphi^h) \nonumber \\
    &= a_T^E(\Pi^\nabla T^h,\Pi^\nabla\varphi^h) + a_T^E(T^h - \Pi^\nabla T^h,\varphi^h - \Pi^\nabla\varphi^h).
\end{align}

The first term represents the energy contribution of the projected functions $\Pi^{\nabla} T^h$ and $\Pi^{\nabla} \varphi^h$, which can be directly computed through polynomial integration. However, projecting variables onto the polynomial space using the projection operator causes energy loss, manifested as rank deficiency in the consistency stiffness matrix, making it non-invertible. The second term represents the energy contribution of projection residuals  $T^h-\Pi^{\nabla} T^h$ and $\varphi^h-\Pi^{\nabla} \varphi^h$. Adding this stabilization term makes the stiffness matrix full rank and positive definite, which is why the second term is called the stabilization term~\cite{RUSSO20161968}.

Therefore, the element temperature stiffness matrix consists of two parts:
\begin{equation}\label{eq:stiffness_decomposition}
    \bm{K}_E^{\text{thermal}} = \bm{K}_E^c + \bm{K}_E^s ,
\end{equation}
where $\bm{K}_E^c$ is the consistency stiffness matrix and $\bm{K}_E^s$ is the stabilization stiffness matrix.
The stabilization stiffness matrix is expressed as~\cite{mengolini2019engineering}
\begin{equation}\label{eq:stabilization_stiffnessT}
    \bm{K}_E^s = \tau^h \text{tr}(\bm{K}_E^c)(\bm{I} - \bm{\Pi}^\nabla)^T(\bm{I} - \bm{\Pi}^\nabla),
\end{equation}
where $\tau^h$ is a user-defined parameter taken as $1/2$ in this work, $\text{tr}(\cdot)$ denotes the trace operator, and $\bm{I}$ is the identity matrix,
$\bm{\Pi}^\nabla$ is the matrix form of the projection operator $\Pi^\nabla$.

\subsection{VEM for coupled thermomechanical analysis}
\label{sec:VEM_thermomechanical}

For elasticity problems, we first define the vector function space for VEM as~\cite{herrera2023numerical}
\begin{equation}\label{eq:vector_space}
    \bm{\mathcal{V}}_{\bm{u}}(E) := \left\{\bm{v}^h \in [\mathcal{H}^1(E)]^2: \bm{v}^h|_{\partial E} \in [C^0(\partial E)]^2, \bm{v}^h|_e \in [\mathcal{M}_k(E)]^2, \Delta\bm{v}^h|_E \in [\mathcal{M}_{k-2}(E)]^2 \right\}.
\end{equation}
where the lowercase $e$ represents the edges on the element boundary $\partial E$, while uppercase $E$ denotes the polygonal element itself. The condition $ \bm{v}^h|_e \in [\mathcal{M}_k(E)]^2$ specifies that the function restriction on each boundary edge $e$ belongs to the polynomial space, while $\Delta\bm{v}^h|_E \in [\mathcal{M}_{k-2}(E)]^2 $defines the Laplacian behavior within the element interior $E$.
With interpolation order $k = 1$, $\mathcal{M}_k(E)$ represents the space of polynomial functions of degree not exceeding $k$ and $\Delta$ represents the Laplace Operator such that $\Delta\boldsymbol{v}^h=\nabla^2\boldsymbol{v}^h=\sum_{i=1}^n\frac{\partial^2\boldsymbol{v}^h}{\partial x_i^2}$.

For elasticity problems, the polynomial space $\bm{\mathcal{M}}_k(E)$ is spanned by six linearly independent vector basis functions (for $k=1$)~\cite{mengolini2019engineering}:
\begin{equation}\label{eq:basis_functions}
    \bm{\mathcal{M}}_k(E) = \left\{
    m_1,
    m_2,
    m_3,
    m_4,
    m_5,
    m_6
    \right\},
\end{equation}
where $m_1=\binom{1}{0}, m_2=\binom{0}{1}, m_3=\binom{-\eta}{\xi}, m_4=\binom{\eta}{\xi}, m_5=\binom{\xi}{0}, m_6=\binom{0}{\eta}$, ${\xi}$ and ${\eta}$ are the normalized coordinates as defined in Eq. \eqref{eq:scaled_coordinates}.

For a polygonal element with $n_v$ vertices, 
we define the local virtual element space $\bm{\mathcal{V}}^h_{\bm{u}}(E)$. 
For two-dimensional vector problems, the number of degrees of freedom for a polygonal element is $n_d = 2n_v$.
Define the projection operator from the function space $\bm{\mathcal{V}}^h_{\bm{u}}(E)$ to the polynomial space $\bm{\mathcal{M}}_k(E)$
\begin{equation}\label{eq:projection_operator}
    \bm{\Pi}_u^\nabla: \bm{\mathcal{V}}_{\bm{u}}^h(E) \to \bm{\mathcal{M}}_k(E),
\end{equation}
the projection operator can be defined based on the orthogonality condition:
\begin{equation}\label{eq:orthogonality_condition}
    a_u^E(\bm{u}^h - \bm{\Pi}_u^\nabla \bm{u}^h, \bm{p}) = 0, \quad \forall \bm{p} \in \bm{\mathcal{M}}_k(E).
\end{equation}

Considering the basic function $\bm{\phi}$ of $\bm{\mathcal{V}}^h_{\bm{u}}(E)$ and $\bm{m}$ of $\bm{\mathcal{M}}_k$,
the orthogonality condition Eq.\eqref{eq:orthogonality_condition} yields
\begin{equation}\label{eq:projection_matrix_relation}
    \bm{M}\tilde{\bm{\Pi}}_u^{*\nabla} = \overline {\bm{B}},
\end{equation}
where $\tilde{\bm{\Pi}}_u^{*\nabla}$ is the Ritz matrix representation of the projection operator $\bm{\Pi}_u^\nabla$~\cite{xu2024high},
and
\begin{equation}
    \bm{M} = \int_E\hat{\bm{\varepsilon}}^T(\bm{m})\hat{\mathbb{D}}\hat{\bm{\varepsilon}}(\bm{m})\ud\Omega,\quad 
    \overline {\bm{B}} = \int_{E}\hat{\bm{\varepsilon}}^T(\bm{m})\hat{\mathbb{D}}\hat{\bm{\varepsilon}}(\bm{\phi})\ud\Omega,
\end{equation}
where $\hat{\square}$ represents the Voigt form and 
\begin{equation}
    \label{r3.varepsilon_m}
    \hat{\bm{\varepsilon}}(\bm{m}) = 
    \begin{bmatrix}
        \frac{\partial m_1}{\partial x} & \cdots & \frac{\partial m_6}{\partial x}\\
        \frac{\partial m_1}{\partial y} & \cdots & \frac{\partial m_6}{\partial y} \\
        \frac{\partial m_1}{\partial x}+\frac{\partial m_1}{\partial y} & \cdots & \frac{\partial m_6}{\partial x}+\frac{\partial m_6}{\partial y}
    \end{bmatrix},\quad k=1
\end{equation}

Besides, the matrix $\overline {\bm{B}}$ can be expressed as:
\begin{equation}\label{eq:Z_matrix_green}
    \overline {\bm{B}} = \int_{E}\hat{\bm{\varepsilon}}^T(\bm{m})\hat{\mathbb{D}}\hat{\bm{\varepsilon}}(\bm{\phi})\ud\Omega 
    =-\int_E\left[\nabla\cdot\left(\hat{\mathbb{D}}\hat{\bm{\varepsilon}}(\bm{m})\right)\right]^T\bm{\phi}\ud\Omega+
        \int_{\partial E}\hat{\bm{\varepsilon}}^T(\bm{m})\hat{\mathbb{D}}\bm{n}\bm{\phi}\ud\Gamma
\end{equation}

Since $\bm{m}_\alpha \in \mathcal{M}_k(E)$, $\nabla \cdot (\mathbb{D}\nabla(\bm{m}_\alpha)) \in \mathcal{M}_{k-2}(E)$,  
we neglect the first term on the right-hand side since $k = 1$. 
Besides, 
\begin{equation}
    \bm{n} =  \begin{bmatrix}
        n_x & 0\\ 0 & n_y\\ n_y & n_x
    \end{bmatrix}.
\end{equation}

Define matrix $\overline {\bm{D}} \in \mathbb{R}^{2n_d \times n_k}$ as the values of polynomial basis $\bm{m}_\alpha$ at the degrees of freedom locations:
\begin{equation}\label{eq:J_matrix}
    \overline{\boldsymbol{D}}=\left[\begin{array}{cccc}
dof_1(m_1) & dof_1(m_2) & \cdots & dof_1(m_6) \\
dof_2(m_1) & dof_2(m_2) & \cdots & dof_2(m_6) \\
\vdots & \vdots & \ddots & \vdots \\
dof_{2n_d}(m_1) & dof_{2n_d}(m_2) & \cdots & dof_{2n_d}(m_6)
\end{array}\right], 
\end{equation}
the matrix $\bm{M}$ can be represented as the product of matrices $\overline {\bm{B}}$ and $\overline {\bm{D}}$:
\begin{equation}\label{eq:M_matrix}
    \bm{M} = \overline {\bm{B}} ~ \overline {\bm{D}}.
\end{equation}

Considering the constraints condition, the projection operator $\bm{\Pi}_u^\nabla$ can be solved lastly based on Eq.~\eqref{eq:projection_matrix_relation}.
For implementation details, including the subroutine named `calculatePi' which demonstrates the process of precisely calculating the $\overline{\boldsymbol{B}}$, $\overline{\boldsymbol{D}}$, and $\bm{M}$ matrices as well as the solution formula for the projection operator, see the MATLAB code provided in the website\footnote{You can find the code from https://www.vemhub.com/code}.
Lastly, the bilinear form within an element can be expressed as: for $\bm{u}^h,\bm{v}^h \in \bm{\mathcal{V}}^h$,
\begin{align}
\label{eq:bilinear_decomposition}
 a_u^E(\bm{u}^h,\bm{v}^h) &= a_u^E(\bm{u}^h - \bm{\Pi}_u^\nabla \bm{u}^h + \bm{\Pi}_u^\nabla \bm{u}^h,\bm{v}^h - \bm{\Pi}_u^\nabla \bm{v}^h + \bm{\Pi}_u^\nabla \bm{v}^h) \nonumber \\
        &= a_u^E(\bm{\Pi}_u^\nabla \bm{u}^h,\bm{\Pi}_u^\nabla \bm{v}^h) + a_u^E(\bm{u}^h - \bm{\Pi}_u^\nabla \bm{u}^h,\bm{v}^h - \bm{\Pi}_u^\nabla \bm{v}^h).
\end{align}

Then the stiffness matrix has the form as 
\begin{equation}\label{eq:stiffness_decomposition_matrix}  
    \bm{K}_E^u = \bm{K}_E^{uc} + \bm{K}_E^{us}.
\end{equation}

The consistency stiffness matrix is expressed as:
\begin{equation}\label{eq:consistency_stiffness}
    \bm{K}_E^{uc} = (\tilde{\bm{\Pi}}_u^{*\nabla})^T\bm{M}\tilde{\bm{\Pi}}_u^{*\nabla}.
\end{equation}

Besides, the stabilization stiffness matrix is expressed as~\cite{mengolini2019engineering}:
\begin{equation}\label{eq:stabilization_stiffness}
    \bm{K}_E^{us} = \tau^h\text{tr}(\bm{K}_E^c)(\bm{I} - \tilde{\bm{\Pi}}_u^{\nabla})^T(\bm{I} - \tilde{\bm{\Pi}}_u^{\nabla}),
\end{equation}
where
\begin{equation}\label{eq:projection_relation}
    \tilde{\bm{\Pi}}_u^{\nabla} = \overline {\bm{D}} \tilde{\bm{\Pi}}_u^{*\nabla}.
\end{equation}
and $\tau^h$ is a user-defined stabilization parameter taken as 1/2 in this work following the established practices in references~\cite{artioli2017arbitrary, mengolini2019engineering}, $\text{tr}(\cdot)$ denotes the trace operator, and $\bm{I}$ is the identity matrix.

Considering the thermal strain, the force vector (see Eq.~\eqref{eq:linear_form_mechanical}) has the form as:
\begin{equation}\label{eq:force_vector}
    \bm{F} = \bm{F}_{\text{ext}} + \bm{F}_{\text{thermal}},
\end{equation}
with $\boldsymbol{F}_{\text{ext}}$ represents the external traction forces applied to the Neumann boundary (with body forces assumed to be zero) and $\boldsymbol{F}_{\text{thermal}}$ represents the thermal load contributions from temperature-induced effects, computed as:
\begin{equation}
    \label{s3.Fthermal}
    \bm{F}_{\text{thermal}} =  (\tilde{\bm{\Pi}}_u^{*\nabla})^T \int_{E} \hat{\bm{\varepsilon}}^T(\bm{m}) \hat{\mathbb{D}} \hat{\bm{\varepsilon}}_{\text{t}}  \, \ud\Omega,
\end{equation}
where $\tilde{\bm{\Pi}}_u^{*\nabla}$ is the Ritz matrix for the projection operator solved from Eq. \eqref{eq:projection_matrix_relation}.
The thermal strain is given in Eq.~\eqref{eq:thermal_strain}, $\hat{\square}$ represents the Voigt form.
Besides, $\hat{\bm{\varepsilon}}(\bm{m})$ is defined in Eq.~\eqref{r3.varepsilon_m}.

\section{Stabilization-free VEM for thermomechanical coupling}
\label{s4}
As shown in Eqs.~\eqref{eq:stabilization_stiffnessT} and \eqref{eq:stabilization_stiffness},
the stabilization term is needed, and some parameters should be introduced to ensure that the stiffness matrix has the correct rank.
The existence of stabilization terms makes the calculation results affected by the user-defined stabilization parameter $\tau^h$.
Therefore, here we present a stabilization-free VEM, which gets rid of the dependence on stabilization terms.

\subsection{Projection in stabilization-free virtual element method}
\label{s4.1}
The basic idea of the stabilization-free virtual element method (SFVEM) is to allow the calculation of the higher order $L_2$ projection for the gradient.
Firstly, we consider a local enlarged enhanced virtual element space
\begin{equation} \label{virtual-space-sfvem}
    \mathcal{V}_{1,l}(E):=\left\{
        v\in\mathcal{H}^1,v|_{\partial E}:v|_e\in\mathcal{P}_1(e),\Delta v\in\mathcal{P}_{l+1}(E)
    \right\}.
\end{equation}

Let $\Pi_{l,E}^0\nabla$ be the $L_2$ projection of the gradient of the function in $\mathcal{V}_{1,l}$ to $\left[\mathcal{P}_l(E)\right]^2$ , where following the notation from reference~\cite{Berrone2021}, the superscript $0$ indicates $L_2$-projection, the subscript $l$ defines the polynomial degree of the gradient projection space required for sufficient stabilization, and $E$ denotes the element. The projection is defined by the orthogonality condition 
\begin{equation}
    \label{p07.s3.1.L2}
    \int_E\bm{p}^T\cdot\Pi_{l,E}^0\nabla v\ud\Omega = \int_E\bm{p}^T\cdot\nabla v\ud\Omega,\quad \bm{p}\in \left[\mathcal{P}_l(E)\right]^2.
\end{equation}
where $\boldsymbol{p}$ is a basis function of the polynomial space $[P_l(E)]^2$, and $P_l(E)$ denotes the space of polynomials of degree $\leq l$ on element $E$.
The right side of Eq.~\eqref{p07.s3.1.L2} can be expanded as 
\begin{equation}
    \label{p07.s3.1.L2.2}
    \int_E\bm{p}^T\cdot\nabla v\ud\Omega=\int_{\partial E}\bm{p}^T\cdot\bm{n}_E v\ud\Gamma-\int_E\left(\text{div}\bm{p}\right)v\ud\Omega.
\end{equation}

Considering the projection operator defined in Eq.~\eqref{eq:projection_orthogonality} and Eq.~\eqref{eq:projection_expanded2-thermal}, and  applying the orthogonality condition with $v$ replacing $\varphi^h$, the last term in Eq.~\eqref{p07.s3.1.L2.2} can be calculated by 
\begin{equation}
\int_E\left(\text{div}\bm{p}\right)v\ud\Omega=\int_E\left(\text{div}\bm{p}\right)\Pi^\nabla v\ud\Omega.
\end{equation}

Besides, the parameter $l$, which defines the polynomial degree required for sufficient stabilization, must satisfy the following condition for the lowest-order case based on the degree-of-freedom constraints established in references~\cite{Berrone2021, Berrone2023}:
\begin{equation}
    \label{s4.l}
    (l+1)(l+2)>n_v-1,
\end{equation}
where $n_v$ is the number of vertices of element $E$. 
According to reference \cite{Xu2023}, the approximated gradient $\nabla v$ and $p$ can be expressed as
\begin{equation}
    \label{s4.p}
    \nabla v=\left(\boldsymbol{N}^p\right)^T \nabla \tilde{\boldsymbol{v}}, \quad \boldsymbol{p}=\left(\boldsymbol{N}^p\right)^T \tilde{\boldsymbol{p}}.
\end{equation}

Then the gradient of the variable $v$ can be approximated as
\begin{equation}
    \label{s4.gradient}
    \nabla v = \Pi_{l,E}^0\nabla v = \left(\bm{N}^p\right)^T\bm{\Pi}^m\tilde{\bm{v}},
\end{equation}
where $N^p=\left[\begin{array}{cc}\boldsymbol{m}_l & 0 \\ 0 & \boldsymbol{m}_l\end{array}\right]=\left[\begin{array}{cccccccccc}1 & \xi & \eta & \cdots & \eta^l & 0 & 0 & 0 & \cdots & 0 \\ 0 & 0 & 0 & \cdots & 0 & 1 & \xi & \eta & \cdots & \eta^l\end{array}\right]$, $\boldsymbol{m}_l$ represents the basis functions of  $l$-th order polynomials, $\boldsymbol{\Pi}^m$ is the matrix form of the $L_2$ projection operator, and $\tilde{\square}$ represents the vector form. If $\tilde{\boldsymbol{v}}$ represents displacement, then $\tilde{\boldsymbol{v}}=\left\{v_{1}^{1}\ v_{1}^{2}\ v_{2}^{1} \ v_{2}^{2}\cdots\ v_{n_{v}}^{1}\ v_{n_{v}}^{2}\right\}$, where $v_i^j$ denotes displacement component $j$ (with $j=1,2$ for $x$, $y$ directions) at node $i$, and $n_v$ is the total number of nodes in the element.

Substituting  Eq.~\eqref{s4.p} and Eq.~\eqref{s4.gradient}  into Eq.~\eqref{p07.s3.1.L2.2}  leads to
\begin{equation}
    \label{p07.s3.1.L2.3}
    \tilde{\boldsymbol{p}}^T \int_E \boldsymbol{N}^p\left(\boldsymbol{N}^p\right)^T \mathrm{~d} \Omega \boldsymbol{\Pi}^m \tilde{\boldsymbol{v}}=\tilde{\boldsymbol{p}} ^T\int_{\partial E}\left(\boldsymbol{N}^p \cdot \boldsymbol{n}\right) \phi^T \mathrm{~d} \Gamma \tilde{\boldsymbol{v}}-\tilde{\boldsymbol{p}} ^T\int_E\left(\operatorname{div} \boldsymbol{N}^p\right) \phi^T \mathrm{~d} \Omega \tilde{\boldsymbol{v}}.
\end{equation}

Since this holds for all $\tilde{\boldsymbol{v}}$ and $\tilde{\boldsymbol{p}}$, Eq.~\eqref{p07.s3.1.L2.3} can be written as
\begin{equation}
    \label{p07.s3.1.L2.4}
\int_E \boldsymbol{N}^p\left(\boldsymbol{N}^p\right)^T \mathrm{~d} \Omega \boldsymbol{\Pi}^m=\int_{\partial E}\left(\boldsymbol{N}^p \cdot \boldsymbol{n}\right) \phi^T \mathrm{~d} \Gamma-\int_E\left(\operatorname{div} \boldsymbol{N}^p\right) \phi^T \mathrm{~d} \Omega.
\end{equation}

If the right side of  Eq.~\eqref{p07.s3.1.L2.4} is computable, then the projection matrix $\boldsymbol{\Pi}^m$ can be calculated by
\begin{equation}
    \label{eq:projection_sfvem}
    \boldsymbol{\Pi}^m=\tilde{\boldsymbol{G}}^{-1} \tilde{\boldsymbol{B}},
\end{equation}
where
\begin{equation}
    \label{eq:G_matrix}
   \tilde{\boldsymbol{G}}=\int_E \boldsymbol{N}^p\left(\boldsymbol{N}^p\right)^T \mathrm{~d} \Omega,
\end{equation}
\begin{equation}
    \label{eq:B_matrix}
   \tilde{\boldsymbol{B}}=\int_{\partial E}\left(\boldsymbol{N}^p \cdot \boldsymbol{n}\right) \phi^T \mathrm{~d} \Gamma-\int_E\left(\operatorname{div} \boldsymbol{N}^p\right) \phi^T \mathrm{~d} \Omega.
\end{equation}
\subsection{Stiffness matrix}
\label{s4.2}
For bilinear form for the thermal conduction problem given in Eq.~\eqref{s2.bilinearT},
the discrete bilinear form without any stabilization term can be easily constructed as 
\begin{equation}
   a_{T}^{E}(\vartheta,w)=\int_{E}\Pi_{l,E}^{0}\nabla\vartheta\cdot\Pi_{l,E}^{0}\nabla w\mathrm{d}\Omega,
\end{equation}
where $\vartheta$ and $w$ represent the temperature field variables in the heat conduction problem.
Considering the gradient in Eq.~\eqref{s4.gradient}, the stiffness matrix for the thermal conduction problem has the form as 
\begin{equation}
    \bm{K}_E^{\text{thermal}} = \lambda_k\left(\bm{\Pi}^m\right)^T\int_E \bm{N}^p \left(\bm{N}^p\right)^T\ud\Omega\bm{\Pi}^m.
\end{equation}

Easy to find that if the appropriate polynomial order $l$ is chosen based on Eq.~\eqref{s4.l}, 
the stiffness matrix does not require any stabilizing terms. Eq.~\eqref{s4.l} ensures the injectivity of the gradient projection operator and guarantees that the projection polynomial space dimension is sufficiently large, avoiding rank deficiency problems in the stiffness matrix~\cite{xu2024stabilization, Berrone2021}.

For the mechanical problem, the strain $\bm{\varepsilon}$ ($\hat{\bm{\varepsilon}}$ in Voigt notation) can be approximated by 
\begin{equation}
    \label{s4.strain}
    \hat{\bm{\varepsilon}}(\bm{u}) = \bm{A}\left[\left(\bm{N}^p\right)^T\otimes\mathbb{I}_2\right]\tilde{\bm{u}},
    = \bm{A}\bm{N}_p^T\bm{\Pi}_m\tilde{\bm{u}},
\end{equation}
where $\otimes$ is the Kronecker product and $\mathbb{I}_2$ represents $2\times 2$ order identity matrix, and
\begin{equation}
    \bm{A} = \begin{bmatrix}
        1 & 0 & 0 & 0\\
        0 & 0 & 0 & 1\\
        0 & 1 & 1 & 0\\
    \end{bmatrix}, \bm{N}_p = \bm{N}^p\otimes\mathbb{I}_2,\bm{\Pi}_m = \bm{\Pi}^m\otimes\mathbb{I}_2.
\end{equation}

Considering the bilinear form in Eq.~\eqref{eq:bilinear_form_mechanical},
the stiffness matrix has the form as 
\begin{equation}
    \bm{K}_E^u = \bm{\Pi}_m^T\int_E \bm{N}_p\bm{A}^T\hat{\mathbb{D}}\bm{A}\bm{N}_p^T\ud\Omega\bm{\Pi}_m.
\end{equation}

Then body force $\bm{F}_{\text{thermal}}$ generated by the thermal strain can be calculated by the same equation as used in Eq.~\eqref{s3.Fthermal}.
Or, considering Eq.~\eqref{s4.strain},
\begin{equation}
    \bm{F}_{\text{thermal}} =  \bm{\Pi}_m^T \int_{E} \bm{N}_p\bm{A}^T \hat{\mathbb{D}} \hat{\bm{\varepsilon}}_{\text{t}}  \, \ud\Omega.
\end{equation}

\section{Mesh generation and optimization strategies} 
\label{sec:meshing}

The VEM offers exceptional mesh flexibility, accommodating elements with any number of nodes, which significantly simplifies coupling between multi-scale components. 
This makes it ideally suited for geometric multiscale electronic packaging models. 
Its ability to handle non-matching mesh structures reduces element count while preserving computational efficiency, convergence, and accuracy. Non-matching meshes refer to the situation where, during the independent discretization of two or more adjacent regions within a computational domain, the positions of interface nodes between these regions do not align, leading to misaligned nodes and edges along their shared interface. 

During model development, polygon meshes can be independently generated for each domain to assembly. 
At component interfaces, each component incorporates nodes from the contact surfaces of adjacent components as new interface nodes. When node coordinates coincide, duplicate nodes are merged to achieve node matching and form new polygonal elements. Alternatively, targeted local mesh refinement can be applied as needed.
This non-matching mesh approach provides enhanced simulation flexibility and prevents local mesh modifications from propagating throughout the global mesh structure.
   \begin{figure}[htbp]
       \centering
       \includegraphics[width=0.5\textwidth]{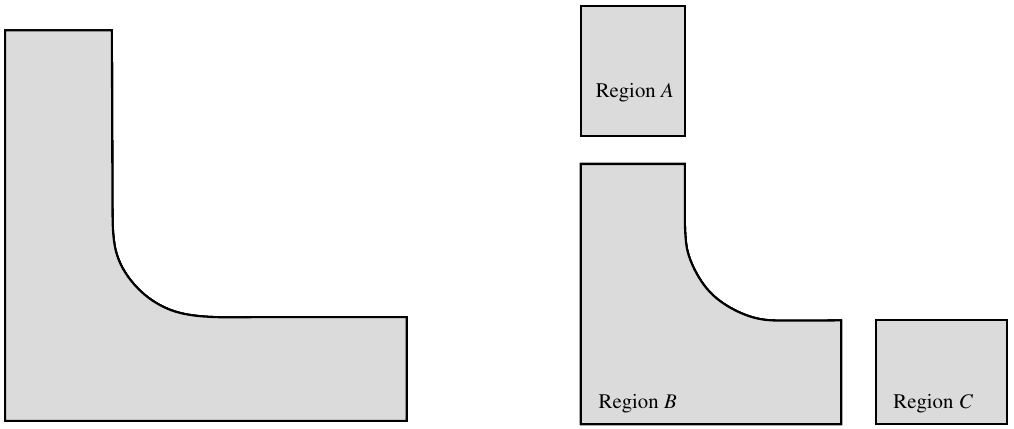}
   \caption{L-shaped structure and domain decomposition strategy.}
   \label{fig:model_partition}
\end{figure}

\begin{figure}[htbp]
   \centering
   \begin{subfigure}[b]{0.46\textwidth}
       \centering
       \includegraphics[width=0.5\textwidth]{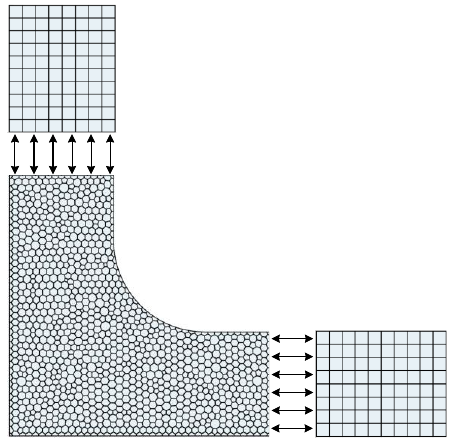}
       \caption{Multi-region mesh distribution}
       \label{fig:mesh_interface_a}
   \end{subfigure}
   \hfill
   \begin{subfigure}[b]{0.52\textwidth}
       \centering
       \includegraphics[width=0.65\textwidth]{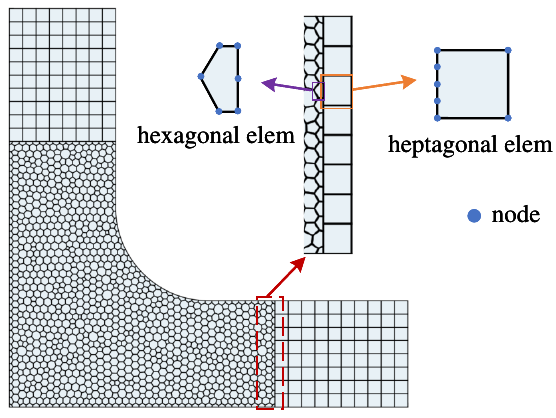}
       \caption{Polygonal-quadrilateral interface detail}
       \label{fig:mesh_interface_b}
   \end{subfigure}
   \caption{Non-matching mesh strategy for L-shaped structure: Region B employs polygonal elements for stress concentration areas, while Regions A and C utilize quadrilateral elements for regular domains.}
   \label{fig:mesh_regions}
\end{figure}

Consider an L-shaped geometric structure, 
as shown in Fig.~\ref{fig:model_partition}. Stress concentration typically occurs at the filleted corner, 
as demonstrated in engineering practice. 
To accurately capture this mechanical behavior, a higher element density is required in the stress concentration region, 
while element density can be reduced elsewhere to enhance computational efficiency. 
Fig.~\ref{fig:model_partition} illustrates our approach where different components are independently modeled and meshed using distinct element types (polygonal and quadrilateral). 
At component interfaces, we adapt non-matching meshes by incorporating interface nodes into polygonal elements. 
As shown in Fig.~\ref{fig:mesh_regions}, the model discretization employs polygonal elements in regions requiring mesh flexibility and quadrilateral elements in regular domains. As illustrated in Fig.~\ref{fig:mesh_interface_b}, the hexagonal element on the left and the quadrilateral element on the right represent different mesh types that can be seamlessly connected through non-matching interface treatment, allowing for flexible mesh transitions without requiring conforming mesh boundaries or transition elements. The interface treatment automatically converts non-matching meshes into polygonal elements based on the spatial distribution of boundary nodes.
 In contrast, FEM relies on standard elements (triangles/quadrilaterals) and requires transition meshes to connect different regions in geometrically multi-scale models, resulting in higher computational costs and mandatory strict node matching at interfaces.

\begin{figure}[htbp]
   \centering
   \includegraphics[width=0.8\textwidth]{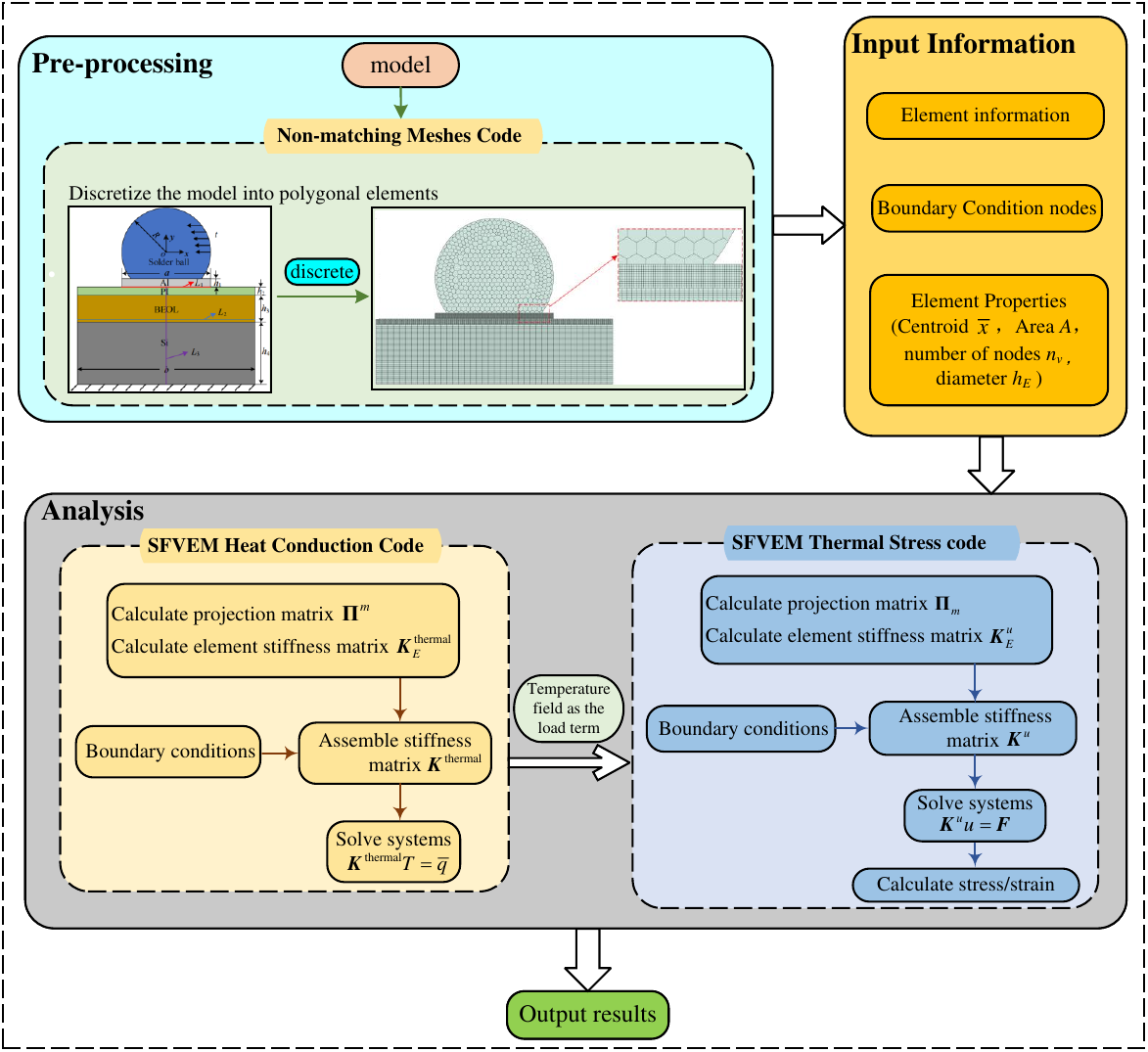}
   \caption{SFVEM-based thermomechanical coupling analysis workflow.}
   \label{fig:vem_flowchart}
\end{figure}

The SFVEM thermomechanical coupling analysis is implemented through two sequential Matlab programs: 
first, a SFVEM heat conduction program calculates the model's temperature field; subsequently, 
a SFVEM thermoelastic program analyzes thermal stresses and displacements using this temperature field as a load input. 
Fig.~\ref{fig:vem_flowchart} illustrates the complete workflow of this thermomechanical coupling analysis scheme based on the virtual element method. 
It should be noted that the sequential nature of our approach (first solving the thermal problem, then applying the temperature field as loading in the mechanical analysis) eliminates the need for iterative convergence criteria typically required in fully coupled nonlinear thermomechanical problems.

\section{Numerical Examples}
\label{Examples}
In this section, we present several numerical examples that demonstrate the effectiveness of our proposed VEM formulation for thermoelastic analysis in electronic packaging. 
These examples include both verification problems with known analytical solutions and practical electronic packaging applications. All numerical examples employ first-order polynomial spaces ( $k=1$ ) for both VEM and SFVEM formulations, as defined in  Eqs.~\eqref{eq:basis_functions} and~\eqref{virtual-space-sfvem}. For SFVEM, the $L_2$ projection order $l$ is determined according to the inequality constraint in Eq.~\eqref{s4.l}; given that the maximum number of element vertices is $12$ across all test cases, we set $l=2$.

Through these selected test cases, we validate the method's accuracy while highlighting its significant advantages in handling complex geometries and non-matching meshes. Sections 6.1 and 6.2 systematically validate the thermomechanical coupling capabilities of SFVEM through comparison with analytical solutions and two-material electronic packaging models, establishing the method's accuracy, convergence, and stability. The subsequent sections (6.3-6.5) focus on mechanical reliability analysis under pure mechanical loading, utilizing SFVEM's mesh flexibility to address geometric multiscale challenges and demonstrate the method's adaptability to complex electronic packaging geometries. It should be noted that interface nodes belong to both materials, with material properties calculated according to their respective elements. Therefore, interface nodes exhibit different stress values depending on the associated material. For validation against finite element results, the stress values at interface nodes are averaged for comparison purposes.

\subsection{Thermoelastic analysis of thick-walled cylindrical structure}
\label{exam1}
\begin{figure}[htbp]
    \centering
    \includegraphics[width=0.9\textwidth]{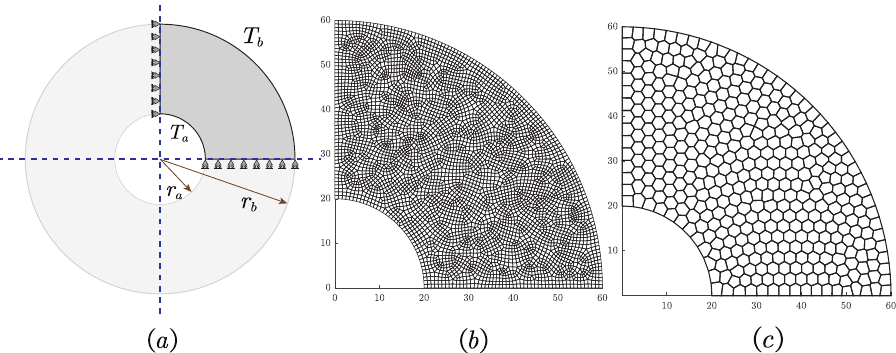}
    \caption{Thick-walled cylinder configuration: (a) Full geometry with inner radius $r_a$, outer radius $r_b$, and applied temperatures $T_a$ and $T_b$; 
    (b) quad mesh, (c) Polygonal mesh.}
    \label{fig:cylinder_model}
\end{figure}

We examine a homogeneous thick-walled cylinder subjected to thermal loading to verify our different VEMs. 
The analytical solution for this model serves as a benchmark for accuracy validation. 
Fig.~\ref{fig:cylinder_model} illustrates the geometry and boundary conditions of this verification case. 
Here, we define a thick-walled cylinder with inner radius $r_a = 20$ mm and outer radius $r_b = 60$ mm. 
The material properties include Young's modulus $E = 460000$ MPa, Poisson's ratio $\nu = 0.3$, 
thermal expansion coefficient $\alpha = 7.4 \times 10^{-6}$ K$^{-1}$, and thermal conductivity $\lambda_k = 20$ W/(m$\cdot$K). 
For the thermal boundary conditions, 
we maintain the inner surface at $T_a = 0$ K and the outer surface at $T_b = 500$ K under steady-state conditions. 
Due to symmetry, we optimize computational efficiency by modeling only one quarter of the cylinder, as illustrated in Fig.~\ref{fig:cylinder_model}.

Given the material parameters and boundary conditions, the analytical temperature field is
\begin{equation}\label{eq:temperature_analytical}
T = T_a + \frac{T_b - T_a}{\ln(r_b/r_a)}\ln\left(\frac{r}{r_a}\right),
\end{equation}
where $r$ represents the radial coordinate. The radial displacement field is given by
\begin{equation}\label{eq:displacement_analytical}
u(r) = B_1r + B_2r^{-1} + Dr\ln r
\end{equation}
The stress components in the thick-walled cylinder can are expressed as
\begin{align}
\sigma_r &= \frac{E}{1-\nu^2}\left(B_1(1+\nu) + B_2(\nu-1)r^{-2} + D\left(\ln r(1+\nu) + 1\right) - \alpha(1+\nu)T(r)\right) \label{eq:radial_stress} \\[8pt]
\sigma_\theta &= \frac{E}{1-\nu^2}\left(B_1(1+\nu) + B_2(1-\nu)r^{-2} + D\left(\ln r(1+\nu) + \nu\right) - \alpha(1+\nu)T(r)\right) \label{eq:hoop_stress}
\end{align}
The thermal coefficient $D$ is
\begin{equation}\label{eq:coefficient_D}
D = \alpha(1+\nu)\frac{T_b - T_a}{2\ln(r_b/r_a)}
\end{equation}
The coefficients $B_1$ and $B_2$ are determined by
\begin{align}
B_1 &= \frac{d_4d_5 - d_2d_6}{d_1d_4 - d_2d_3} \label{eq:coefficient_B1} \\[8pt]
B_2 &= \frac{d_1d_6 - d_3d_5}{d_1d_4 - d_2d_3} \label{eq:coefficient_B2}
\end{align}
where the auxiliary coefficients are defined as
\begin{align} 
d_1 &= 1 + \nu \label{eq:d1} \\[4pt]
d_2 &= r_a^{-2}(\nu - 1) \label{eq:d2} \\[4pt]
d_3 &= 1 + \nu \label{eq:d3} \\[4pt]
d_4 &= r_b^{-2}(\nu - 1) \label{eq:d4} \\[4pt]
d_5 &= -D(\ln(r_a) + 1 + \nu\ln(r_a)) + \alpha(1 + \nu)T_a \label{eq:d5} \\[4pt]
d_6 &= -D(\ln(r_b) + 1 + \nu\ln(r_b)) + \alpha(1 + \nu)T_b \label{eq:d6}
\end{align}

\begin{figure}[htbp]
    \centering
    \includegraphics[width=0.6\textwidth]{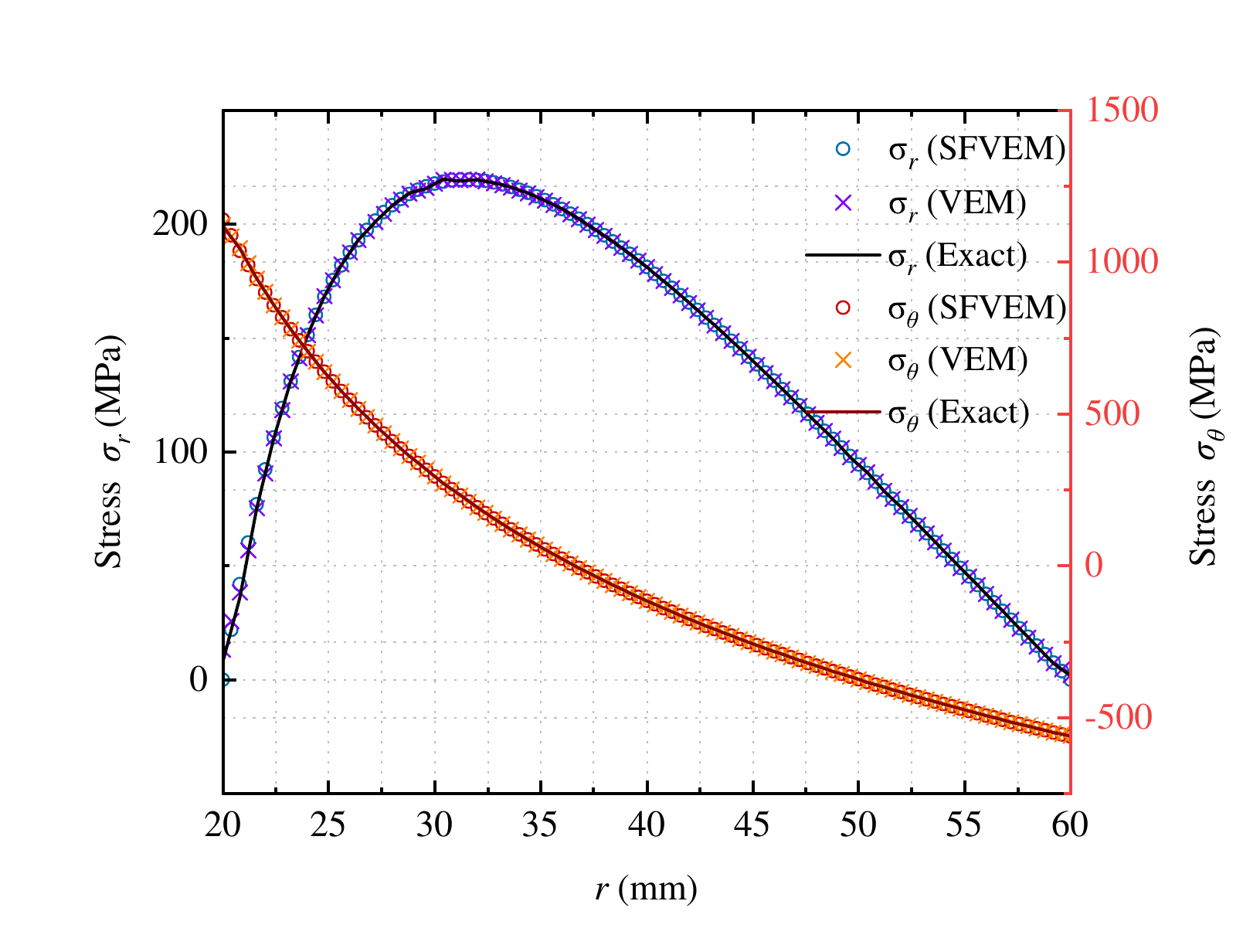}
    \caption{Radial ($\sigma_r$) and circumferential ($\sigma_\theta$) stress distributions under thermal load.}
    \label{fig:EX1_stress_distribution}
\end{figure}

Using the prescribed temperature boundary conditions, 
we first compute the temperature field and then apply it as a load to calculate the displacement and stress distribution. 
The mesh used here is the quad mesh which contains 5,073 nodes illustrated in Fig.~\ref{fig:cylinder_model}.
The calculated radial stress $\sigma_r$ and circumferential stress $\sigma_\theta$ are given in Fig.~\ref{fig:EX1_stress_distribution}.
The results obtained by VEM and SFVEM all have excellent agreement with the analytical solution.

To quantitatively assess the computational accuracy, 
we compare our VEM results with finite element analysis using identical mesh density. 
We evaluate performance using the average error $E_{\text{AV}}$, defined as
\begin{equation}\label{eq:average_error}
\left\{
    \begin{aligned}
        &E_{\text{AV}}^r = \frac{1}{N}\sum\limits_{n=1}^N \left|\frac{\sigma_r - \sigma_r^{\text{exact}}}{\sigma_r^{\text{exact}}}\right| \times 100 \\
        &E_{\text{AV}}^\theta = \frac{1}{N}\sum\limits_{n=1}^N \left|\frac{\sigma_\theta - \sigma_\theta^{\text{exact}}}{\sigma_\theta^{\text{exact}}}\right| \times 100
    \end{aligned}
\right.
\end{equation}
where $E_{\text{AV}}^r$ and $E_{\text{AV}}^\theta$ represent the average percentage errors in radial and circumferential stresses, respectively.

Tab.~\ref{tab:error_comparison} summarizes the average errors for both numerical methods (quad mesh for VEM and FEM).
The results demonstrate that our VEM approach achieves comparable accuracy to the established finite element method.

\begin{table}[htbp]
    \small
    \centering
    \caption{Comparison of average percentage errors between FEM and VEM methods.}
    \begin{tabular}{lcccc}
        \toprule
        \multirow{2}{*}{Error} & \multicolumn{2}{c}{FEM} & \multicolumn{2}{c}{VEM} \\
        \cmidrule(lr){2-3} \cmidrule(lr){4-5}
        & $E_{\text{AV}}^r$ & $E_{\text{AV}}^\theta$ & $E_{\text{AV}}^r$ & $E_{\text{AV}}^\theta$ \\
        \midrule
            & 0.708\% & 0.36158\% & 0.708\% & 0.36157\% \\
        \bottomrule
    \end{tabular}
    \label{tab:error_comparison}
\end{table}

We next assess the method's convergence characteristics using the root mean square (RMS) $L^2$ error as our accuracy metric. 
This normalized error measure is defined as
\begin{equation}\label{eq:rms_error}
\varepsilon_{\text{rms}} = \frac{1}{\max(|X_e|)}\sqrt{\frac{1}{N}\sum_{i=1}^N |X - X_e|^2}
\end{equation}
where $N$ represents the number of considered points, while $X_e$ and $X$ denote the analytical and numerical results, respectively.
\begin{figure}[htbp]
    \centering
    \begin{subfigure}[b]{0.48\textwidth}
        \centering
        \includegraphics[width=\textwidth]{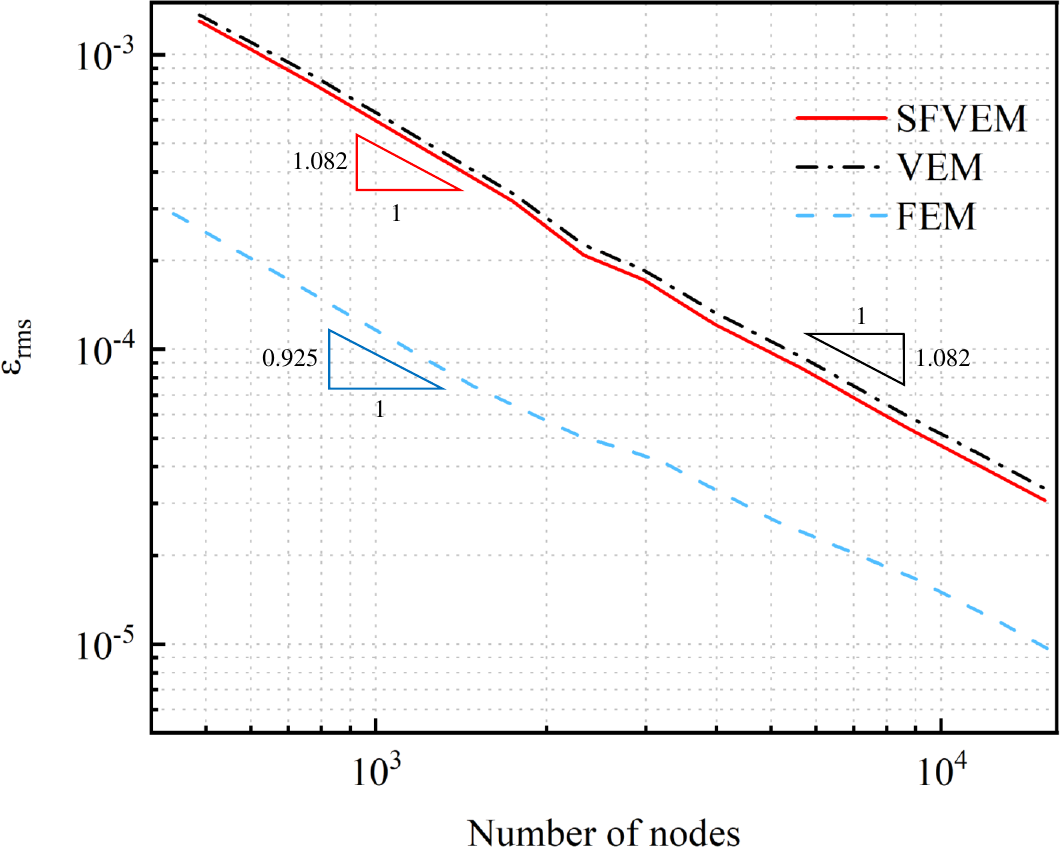}
        \caption{Temperature accuracy versus mesh refinement}
        \label{fig:EX1_temperature_convergence}
    \end{subfigure}
    \hfill
    \begin{subfigure}[b]{0.48\textwidth}
        \centering
        \includegraphics[width=\textwidth]{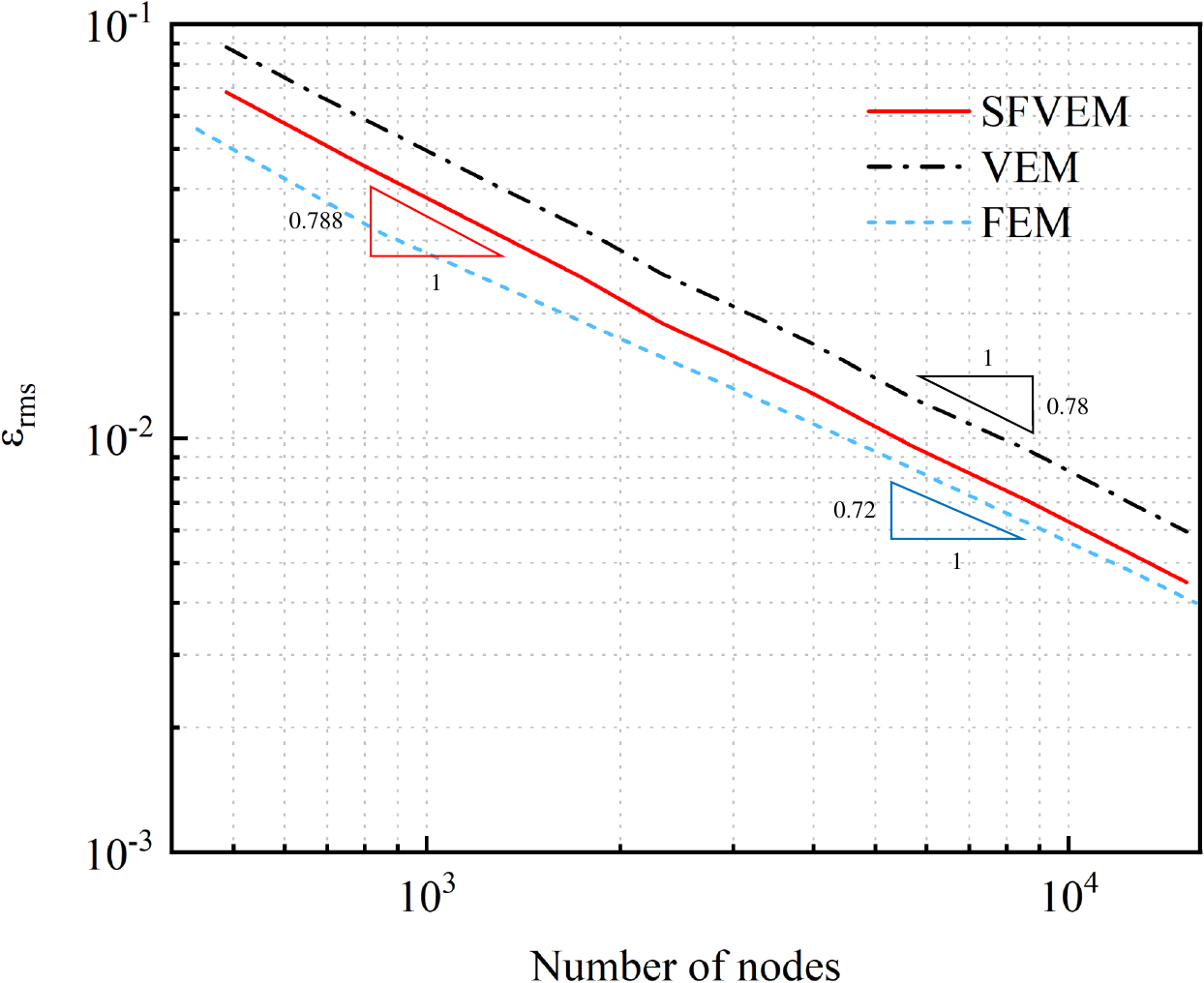}
        \caption{Stress accuracy versus mesh refinement}
        \label{fig:EX1_stress_convergence}
    \end{subfigure}
    \caption{RMS errors with increasing degrees of freedom for (a) temperature field and (b) stress field.}
    \label{fig:EX1_convergence_study}
\end{figure}

Figs.~\ref{fig:EX1_temperature_convergence} and \ref{fig:EX1_stress_convergence}    present the convergence analysis comparing SFVEM, traditional VEM, and conventional FEM, showing RMS errors for temperature and stress fields as functions of degrees of freedom. The convergence rate analysis demonstrates that SFVEM achieves slopes of 1.082 and 0.788 for temperature and stress accuracy respectively, which are comparable to conventional FEM while providing enhanced accuracy compared to traditional VEM. Both figures confirm the robust convergence characteristics and competitive performance of our proposed SFVEM algorithm for thermoelastic problems.

\begin{figure}[htbp]
    \centering
    \begin{subfigure}[b]{0.48\textwidth}
        \centering
        \includegraphics[width=0.8\textwidth]{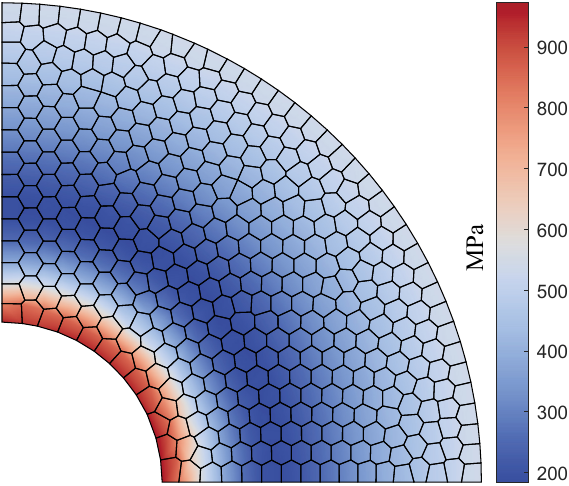}
       \caption{von Mises stress}
        \label{fig:EX1_von Mises stress}
    \end{subfigure}
    \hfill
    \begin{subfigure}[b]{0.48\textwidth}
        \centering
        \includegraphics[width=0.8\textwidth]{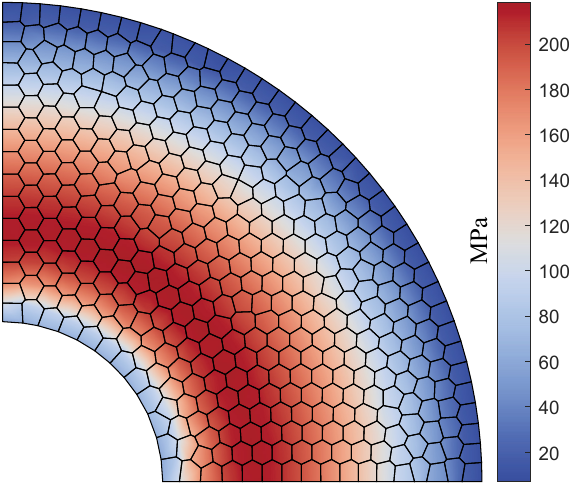}
        \caption{ radial stress $\sigma_r$}
        \label{fig:EX1_radial stress}
    \end{subfigure}
    \caption{Stresses under thermal load:(a) von-Mises stress, (b) radial stress $\sigma_r$ }
    \label{fig:ex1_stress}
\end{figure}

Finally, using the polygonal mesh configuration, the stress contour plots are presented in Fig.~\ref{fig:ex1_stress}. The smooth and accurate stress contours demonstrate the VEM’s excellent suitability for polygonal mesh discretizations.

\subsection{Through-Silicon Via with Copper filling: mechanical and thermal analysis}
\label{exam2}
\begin{figure}[htbp]
    \centering
    \includegraphics[width=0.8\textwidth]{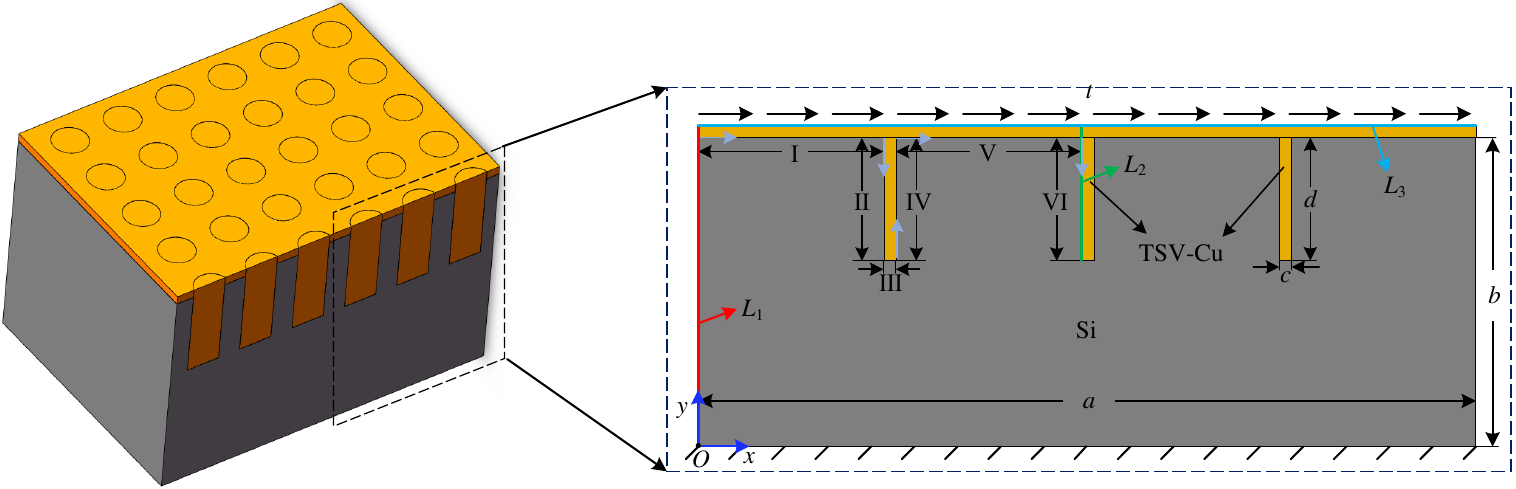}
    \caption{Schematic representation of TSV-Cu structural configuration and dimensions.}
    \label{fig:EX2_tsv_model}
\end{figure}

In this example, we apply our  SFVEM algorithms to analyze both linear elastic and thermoelastic behaviors of a TSV-Cu (Through-Silicon Via with Copper filling) structure.
This critical electronic packaging technology creates blind holes in silicon wafers using dry etching processes, 
then fills them with electroplated copper to provide both mechanical support and electrical interconnection for transistors. 
Fig.~\ref{fig:EX2_tsv_model} illustrates our simplified TSV-Cu structural model. 
In this model, the TSV-Cu components are uniformly distributed along the length direction, with a spacing of $150\, \mu\mathrm{m}$ between each TSV. Each component features a width of $c = 10\,  \mu\mathrm{m}$ and a depth of $d = 100\, \mu\mathrm{m}$, complemented by a top surface coating with a thickness of $10 \,\mu\mathrm{m}$.
The silicon substrate has dimensions of $a = 630\, \mu\text{m}$ and $b = 250\,\mu\text{m}$. 
Table \ref{tab:EX1_material_properties} summarizes the material properties used in our simulation.

\begin{table}[htbp]
    \small
    \centering
    \caption{Mechanical and thermal properties of materials in TSV-Cu structure.}
    \label{tab:EX1_material_properties}
    \begin{tabular}{lcccc}
        \toprule
        Material & Young's modulus & Poisson's ratio & Thermal expansion & Conductivity \\
        & $E$(MPa) & $\nu$ & $\alpha$(K$^{-1}$) & $\lambda_k$(W/m·K) \\
        \midrule
        Cu & 155000 & 0.3 & 17$\times$10$^{-6}$ & 397 \\
        Si & 140000 & 0.25 & 2.8$\times$10$^{-6}$ & 149 \\
        \bottomrule
    \end{tabular}
\end{table}

\subsubsection{Linear elastic analysis of TSV-Cu structure}
We begin by examining the linear elastic response of the TSV-Cu structure using our SFVEM. 
As shown in Fig.~\ref{fig:EX2_tsv_model}, we apply fixed constraints ($u_x = u_y = 0$) to the model's bottom surface, while imposing a shear load of  $7$ N/mm on the upper surface of the copper coating layer. Additionally, the left and right boundaries of the model are set as free boundaries.

\begin{figure}[htbp]
    \centering
    \begin{subfigure}[b]{0.48\textwidth}
        \centering
        \includegraphics[width=0.8\textwidth]{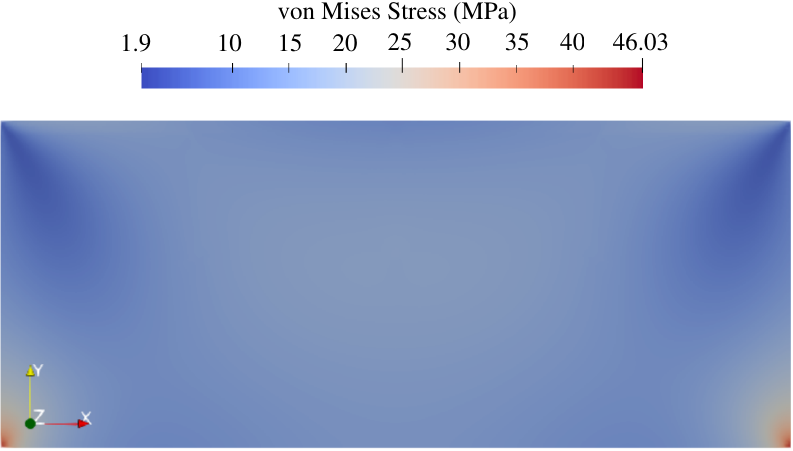}
       \caption{SFVEM-computed stress field}
        \label{fig:EX2_von_mises_vem}
    \end{subfigure}
    \hfill
    \begin{subfigure}[b]{0.48\textwidth}
        \centering
        \includegraphics[width=0.8\textwidth]{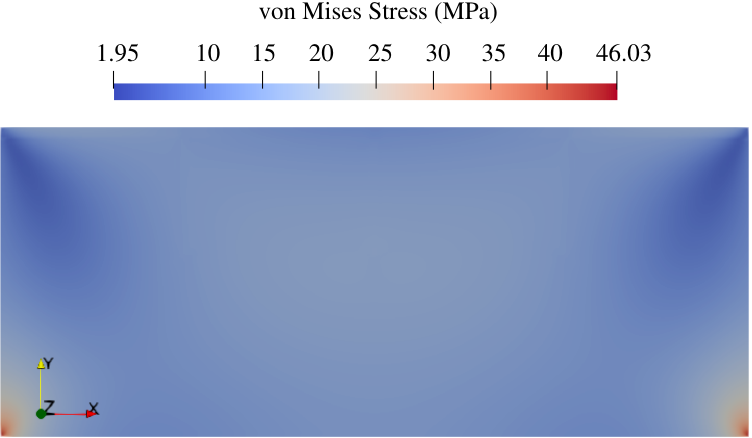}
        \caption{FEM-computed stress field}
        \label{fig:EX2_von_mises_fem}
    \end{subfigure}
    \caption{von Mises stress distribution comparison between (a) SFVEM (polygonal mesh) and (b) FEM (triangular mesh).}
    \label{fig:EX2_stress_comparison}
\end{figure}

\begin{figure}[htbp]
    \centering
    \begin{subfigure}[b]{0.48\textwidth}
        \centering
        
        \includegraphics[width=0.8\textwidth]{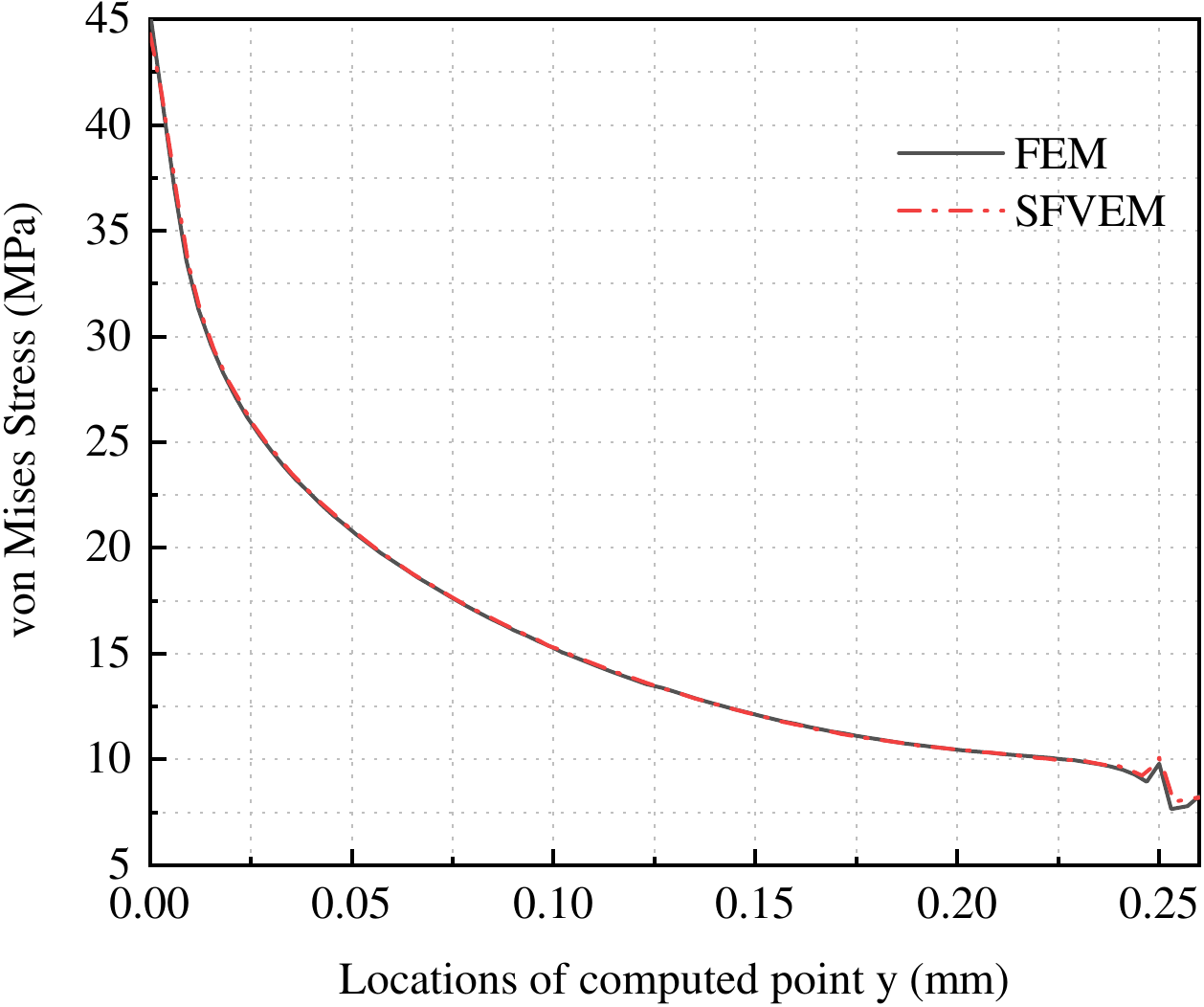}
        \caption{Results along the line $L_1$}
        \label{fig:EX2_von_mises_tsvl1}
    \end{subfigure}
    \hfill
    \begin{subfigure}[b]{0.48\textwidth}
        \centering
        \includegraphics[width=0.82\textwidth]{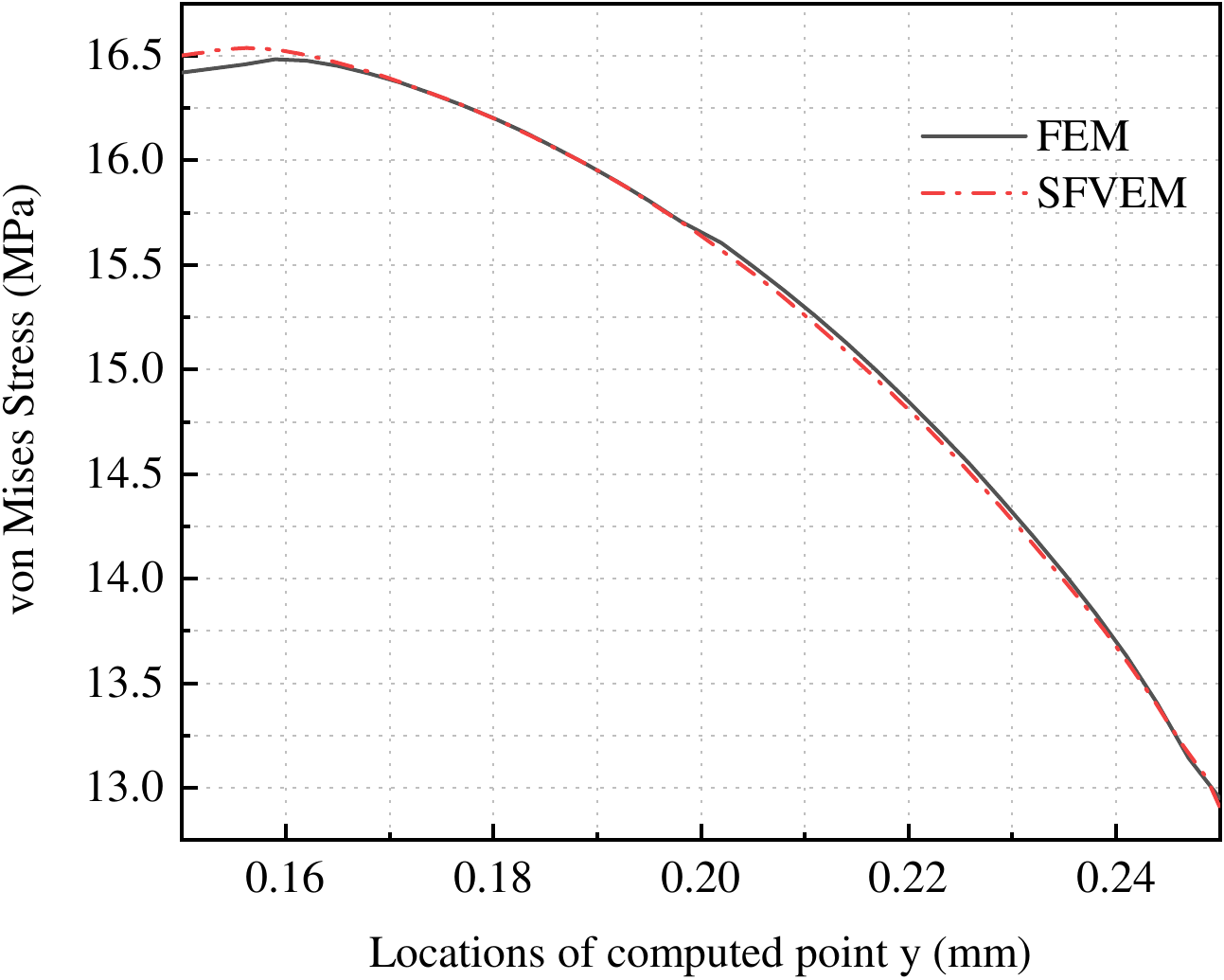}
        \caption{Results along the line $L_2$}
        \label{fig:EX2_von_mises_tsvl2}
    \end{subfigure}
    \hfill
    \begin{subfigure}[b]{0.48\textwidth}
        \centering
        \includegraphics[width=0.8\textwidth]{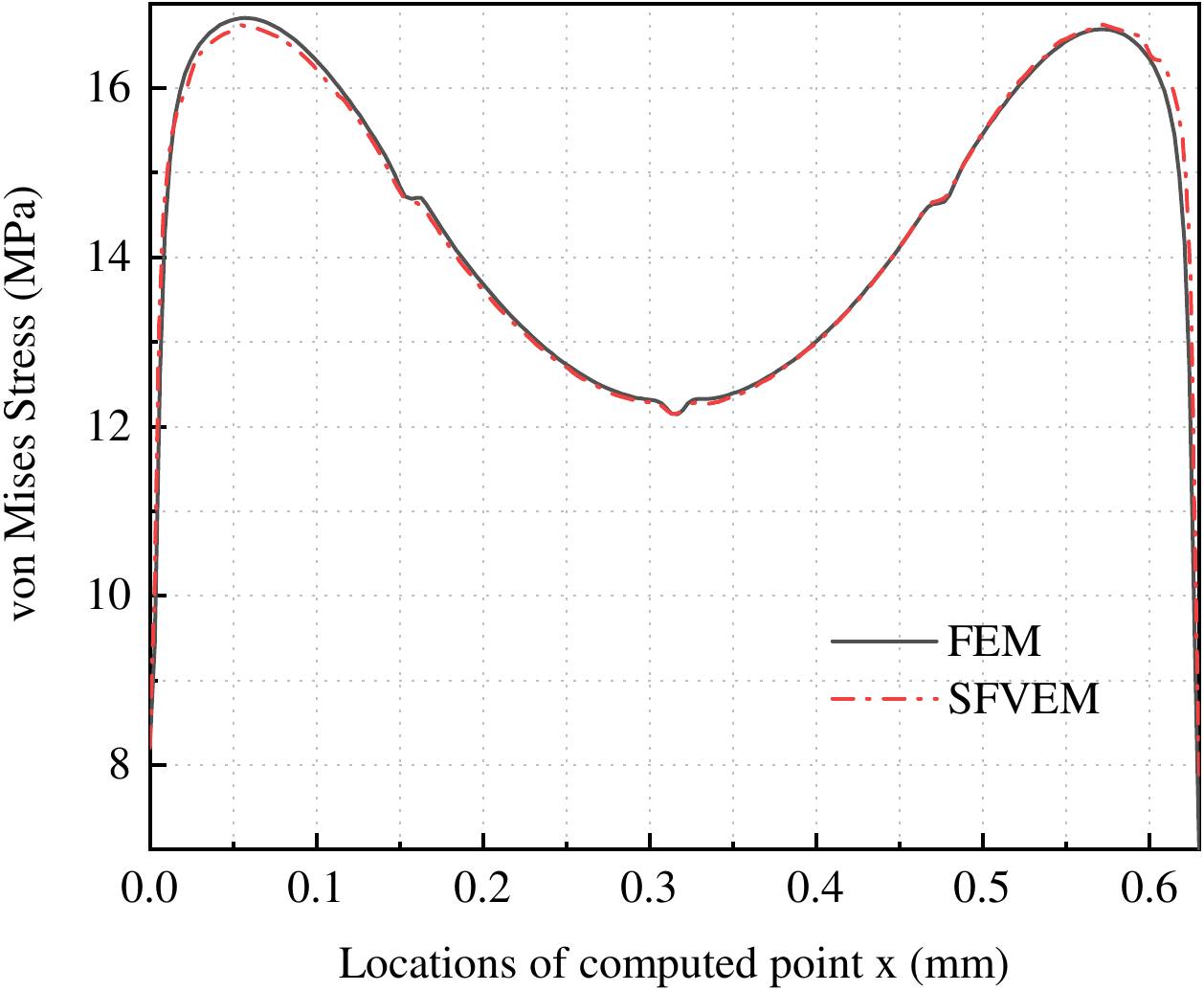}
        \caption{Results along the line $L_3$}
        \label{fig:EX2_von_mises_tsvl3}
    \end{subfigure}
    \caption{Stress distribution comparisons along the lines $L_1$, $L_2$ and $L_3$ in TSV-Cu model.}
    \label{fig:EX2_stress_comparison_tsvl}
\end{figure}

To validate our proposed algorithm, we compare its results with those from commercial software Abaqus, 
which serves as our reference solution. It should be noted that while Fig.~\ref{fig:EX2_tsv_model} shows a periodic structure, the present analysis considers only a finite portion of the 2D cross-section, without implementing periodic boundary conditions. This modeling choice leads to the stress concentrations observed at the left and right boundaries.
Figs.~\ref{fig:EX2_von_mises_vem} and~\ref{fig:EX2_von_mises_fem} present the stress distribution contours from both methods, 
demonstrating excellent agreement between our SFVEM solution and the reference. To further analyze the computational accuracy, three lines ($L_1$, $L_2$, $L_3$) are selected in Fig.~\ref{fig:EX2_tsv_model}, and the stress values along these paths are extracted from the stress distribution contours shown in  Fig.~\ref{fig:EX2_stress_comparison} and compared between FEM and SFVEM. The comparison results are presented in Figs.~\ref{fig:EX2_von_mises_tsvl1}--\ref{fig:EX2_von_mises_tsvl3}.


\begin{figure}[htbp]
    \centering
    \includegraphics[width=0.45\textwidth]{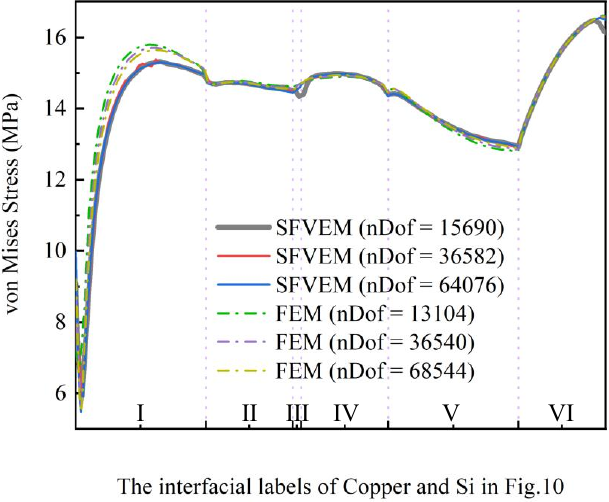}
    \caption{von Mises stress variation along the TSV-Cu/Si interface at different mesh refinement levels.}
    \label{fig:EX2_interface_stress}
\end{figure}

\begin{figure}[htbp]
    \centering
    \begin{subfigure}[b]{0.9\textwidth}
        \centering
        \includegraphics[width=1\textwidth]{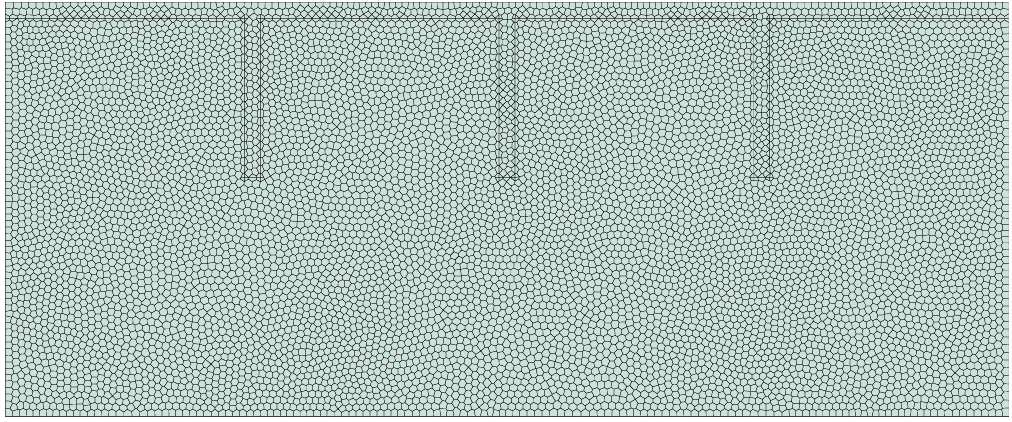}
       \caption{Polygonal mesh discretization of the TSV model }
        \label{fig:EX2_SFVEM MESH}
    \end{subfigure}
    \hfill
    \begin{subfigure}[b]{0.9\textwidth}
        \centering
        \includegraphics[width=1\textwidth]{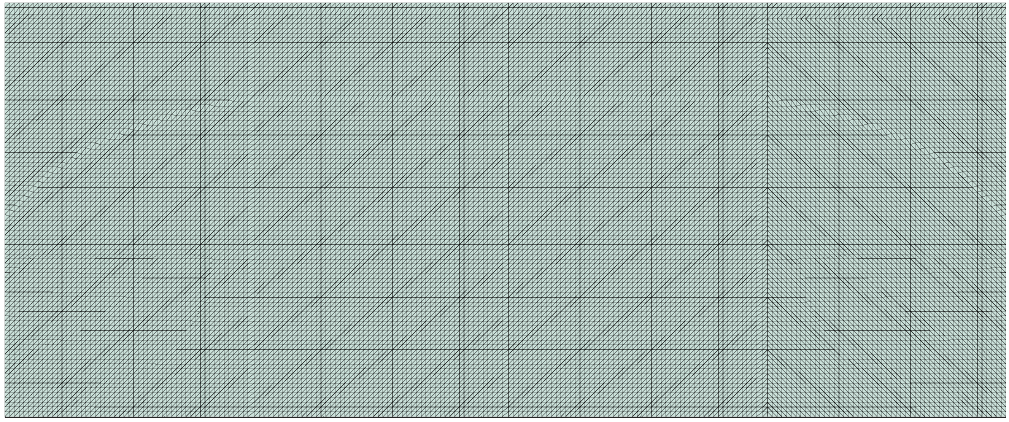}
        \caption{Triangular mesh discretization of the TSV model for FEM }
        \label{fig:EX2_FEM MESH}
    \end{subfigure}
    \caption{The mesh discretization of TSV-Cu structure.}
    \label{fig:EX2_mesh}
\end{figure}

Fig.~\ref{fig:EX2_SFVEM MESH} shows the polygonal mesh discretization of the TSV model used in this analysis. For two-dimensional problems, each vertex possesses two displacement degrees of freedom (in the $x$ and $y$ directions), so the total number of degrees of freedom (nDof) referenced in the following convergence study is twice the number of vertices.
Fig.~\ref{fig:EX2_interface_stress} illustrates the von Mises stress distribution along the TSV-Cu/Si interface across six distinct regions (I through VI), 
calculated using SFVEM at various mesh refinement levels. 
The comparison confirms strong agreement between current results and the FEM solutions. 
However, a larger difference between FEM and SFVEM is observed in region I, which can be attributed to the number of  degrees of freedom in the $y$-direction of the TSV of FEM is less than that of SFVEM, resulting in reduced accuracy of the finite element solution. This discrepancy is amplified in region I due to the significant stress variations present in this area. As evident from the convergence plots, the FEM results progressively approach the SFVEM solution as the degrees of freedom increase.
The smooth stress transition along the interface highlights the proposed algorithm's effectiveness in accurately capturing interfacial stress distributions.


\subsubsection{Thermoelastic response analysis of TSV-Cu structure}  
\begin{figure}[htbp]
    \centering
    \includegraphics[width=0.7\textwidth]{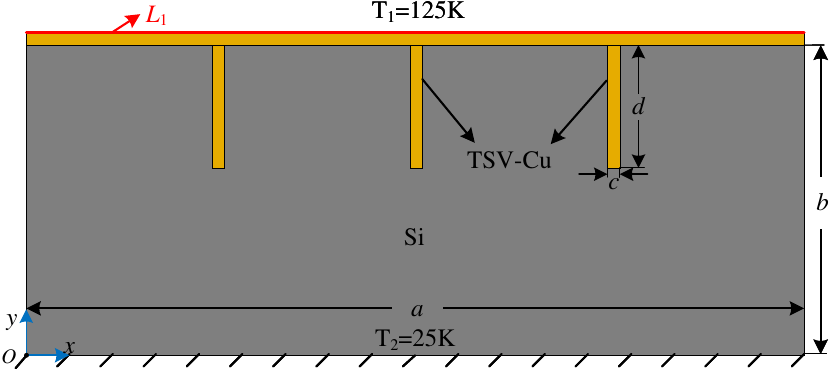}
    \caption{Thermal boundary conditions for TSV-Cu analysis: 125 K at upper surface, 25 K at lower surface, and thermally insulated lateral boundaries.}
    \label{fig:EX2_thermal_boundary}
\end{figure}

Here, the thermoelastic response of the TSV-Cu structure is examined using the proposed SFVEM algorithm. 
Utilizing the same structural model shown in Fig.~\ref{fig:EX2_tsv_model}, 
we apply thermal boundary conditions as illustrated in Fig.~\ref{fig:EX2_thermal_boundary}. 
We maintain fixed displacement constraints ($u_x = u_y = 0$) at the model's bottom surface while imposing a temperature of 125 K on the upper surface and 25 K on the lower surface. 
The left and right boundaries are thermally insulated with zero heat flux ($\partial T/ \partial n = 0$) and are characterized by traction-free boundary conditions.

\begin{figure}[htbp]
    \centering
    \includegraphics[width=0.45\textwidth]{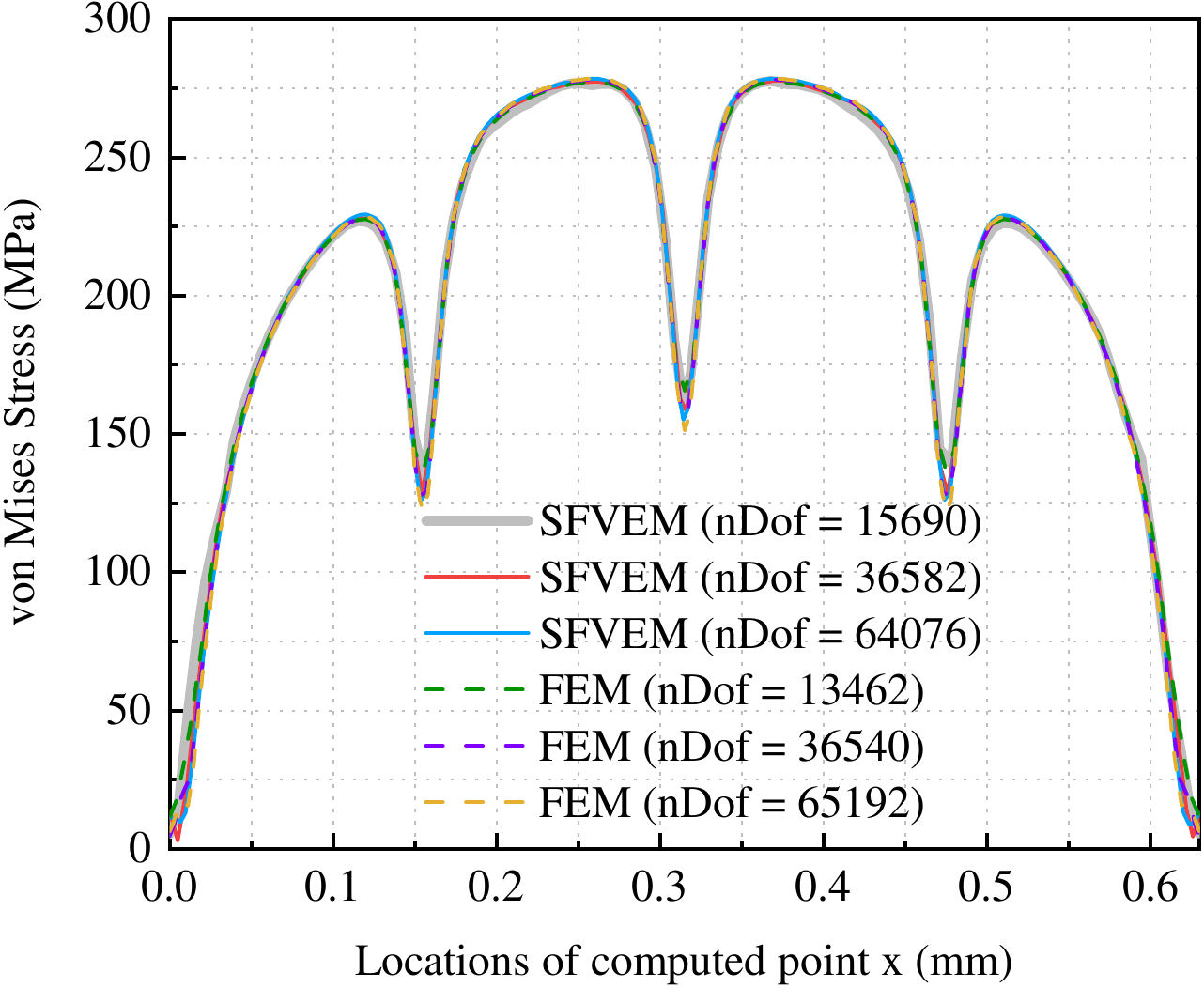}
    \caption{Comparative analysis of thermal stress distribution along the upper surface of TSV-Cu structure at different mesh refinement levels.}
    \label{fig:EX2_thermal_stress}
\end{figure}

To validate the SFVEM implementation for thermoelastic problems, 
we benchmark against commercial software Abaqus. Fig.~\ref{fig:EX2_thermal_stress} compares the thermal stress distributions along the TSV-Cu structure's upper surface (the line $L_1$) calculated by both methods at various mesh refinement levels. 
The results demonstrate excellent agreement between our SFVEM solution and the reference FEM solution across all degrees of freedom.

\subsection{Mechanical analysis of Ball Grid Array structure}
\begin{figure}[htbp]
    \centering
    \includegraphics[width=0.4\textwidth]{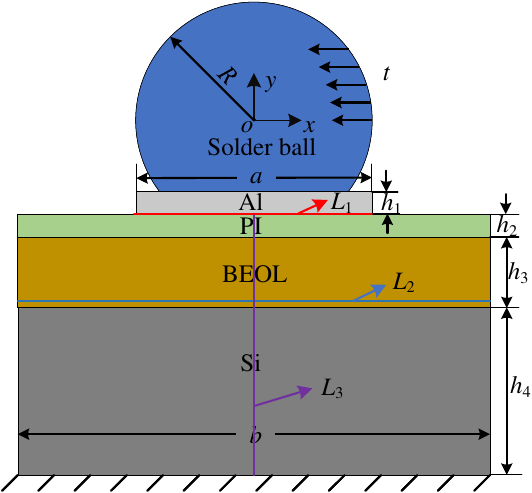}
    \caption{Schematic diagram of BGA model showing solder ball array and package structure.}
    \label{fig:EX3_bga_model}
\end{figure}

Ball Grid Array (BGA) packaging technology represents an advanced high-density surface mount approach utilizing arrays of solder ball interconnections. 
BGA packages deliver two to three times greater memory capacity without volume increase, while providing significant advantages in form factor reduction, thermal management, and electrical performance. 
The BGA model illustrated in Fig.~\ref{fig:EX3_bga_model} has the following geometric parameters: $a = 100\, \mu \text{m}$, $b = 200\,\mu \text{m}$, ${h_1} = {h_2} = 5\,\mu\text{m}$, ${h_3} = 15\,\mu\text{m}$, ${h_4} = 35\,\mu\text{m}$, and $R = 50\,\mu\text{m}$. The contact angle of the Solder ball with the Al plate is $53^{\circ}$.
Tab.~\ref{tab:EX3_bga_properties} summarizes the material properties used in our simulation. 
For boundary conditions, as shown in Fig.~\ref{fig:EX3_bga_model}, we apply fixed constraints at the model's bottom surface while imposing a stress of $t = 10$ MPa on the right side of the solder ball. The vertical coordinate range of the applied load $t$ is applied is $(-0.005, 0.031)~\text{mm}$. All other boundaries are traction-free.

\begin{table}[htbp]
    \small
    \centering
    \caption{Material properties of BGA model components.}
    \label{tab:EX3_bga_properties}
    \begin{tabular}{lcc}
        \toprule
        Material & Young's modulus (MPa) & Poisson's ratio \\
        \midrule
        Solder ball & 10000 & 0.4 \\
        Al & 70000 & 0.3 \\
        PI & 4000 & 0.2 \\
        BEOL & 90000 & 0.3 \\
        Si & 130000 & 0.28 \\
        \bottomrule
    \end{tabular}
\end{table}

\begin{figure}[htbp]
    \centering
    \begin{subfigure}[b]{\textwidth}
        \centering
    \includegraphics[width=0.9\textwidth]{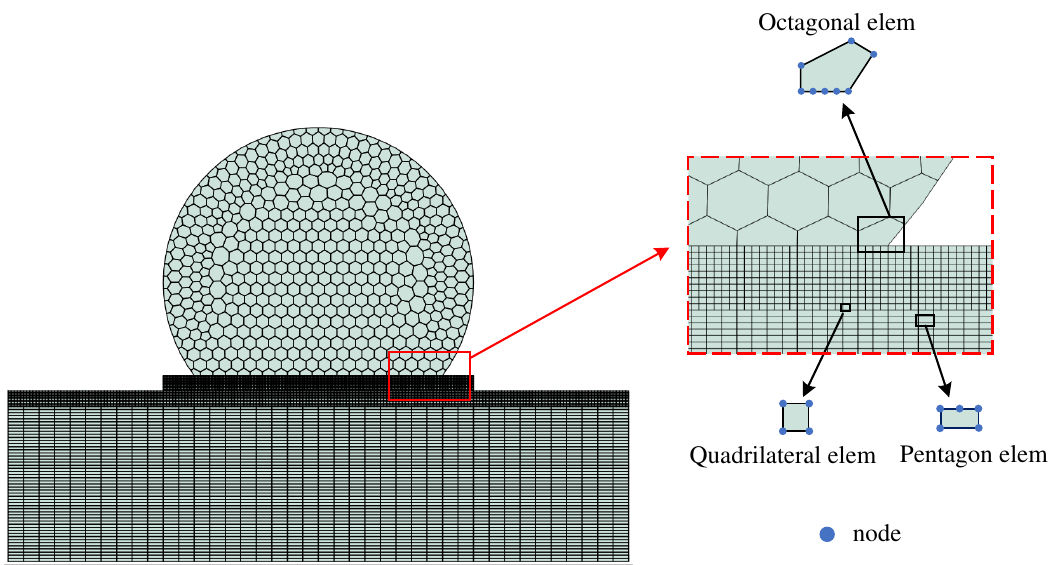}
    \caption{Polygonal mesh discretization of BGA structure for SFVEM.}
    \label{fig:EX3_bga_mesh}
 \end{subfigure}
    \hfill
    \begin{subfigure}[b]{0.85\textwidth}
        \centering
    \includegraphics[width=0.9\textwidth]{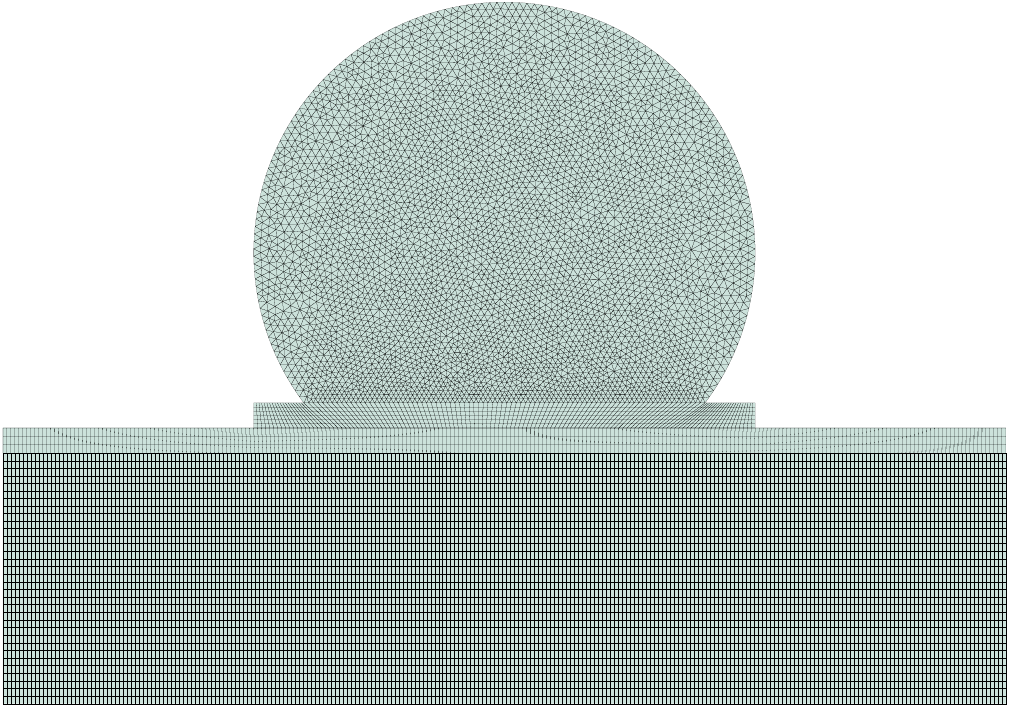}
    \caption{Quad and triangle hybrid mesh discretization of BGA structure for FEM.}
    \label{fig:EX3_bga_mesh_fem}
\end{subfigure}
    \caption{ The mesh discretization of BGA structure.}
\end{figure}

The exceptional mesh flexibility of SFVEM, which accommodates polygonal elements with any number of nodes, enables us to effectively address the complex meshing challenges presented by this BGA model. 
With significant geometric variations between components and anticipated stress concentration at the corner interface between the solder ball and Al layer, 
our approach strategically increases element density near these critical interfaces while reducing mesh density in less critical regions to optimize computational efficiency. 
Fig.~\ref{fig:EX3_bga_mesh} illustrates the polygonal mesh discretization for the BGA model. 
By employing component-specific element sizing and non-matching mesh interfaces, we eliminate the need for transition elements, significantly enhancing computational efficiency.

\begin{figure}[htbp]
    \centering
    \begin{subfigure}[b]{0.48\textwidth}
        \centering
        \includegraphics[width=0.8\textwidth]{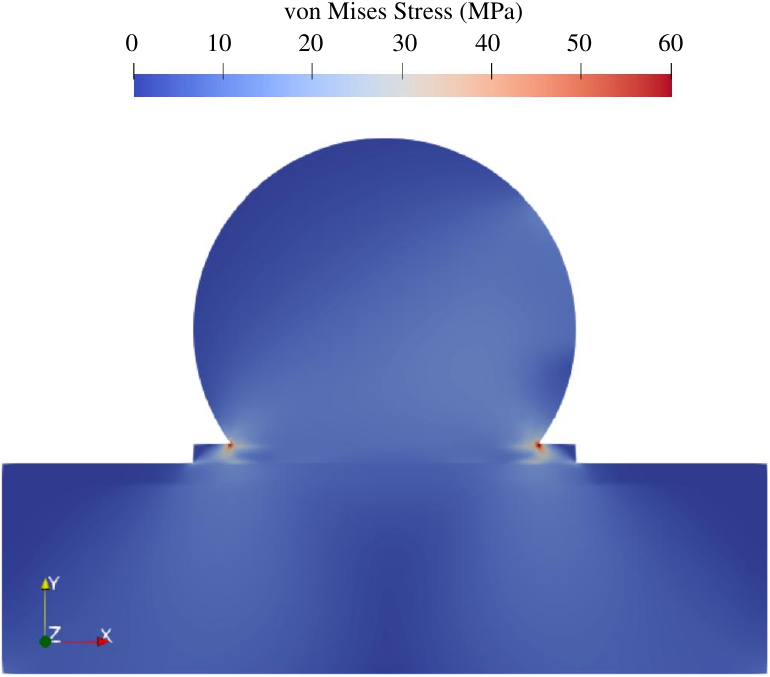}
        \caption{von Mises stress obtained by SFVEM}
        \label{fig:EX3_von_mises_vem_bga}
    \end{subfigure}
    \hfill
    \begin{subfigure}[b]{0.48\textwidth}
        \centering
        \includegraphics[width=0.8\textwidth]{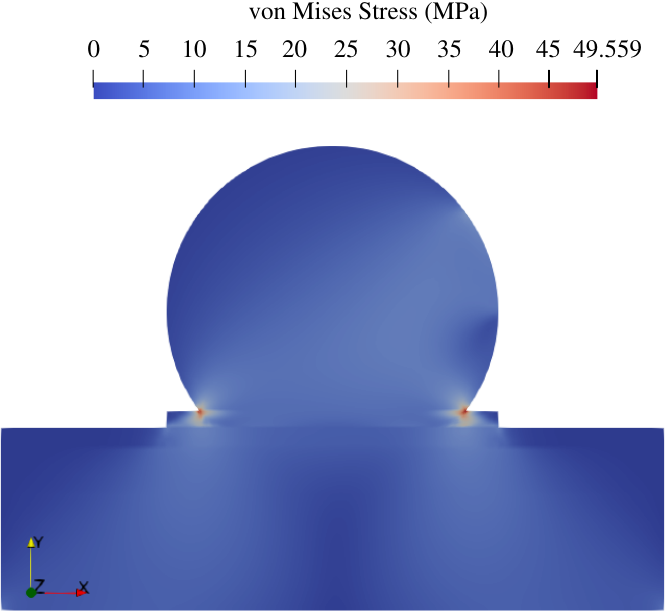}
        \caption{von Mises stress obtained by FEM}
        \label{fig:EX3_von_mises_fem_bga}
    \end{subfigure}
    \caption{Comparison of stress distributions in BGA structure: (a) SFVEM result; (b) FEM result.}
    \label{fig:EX3_stress_comparison_bga1}
\end{figure}

\begin{figure}[htbp]
    \centering
    \begin{subfigure}[b]{0.48\textwidth}
        \centering
        \includegraphics[width=0.85\textwidth]{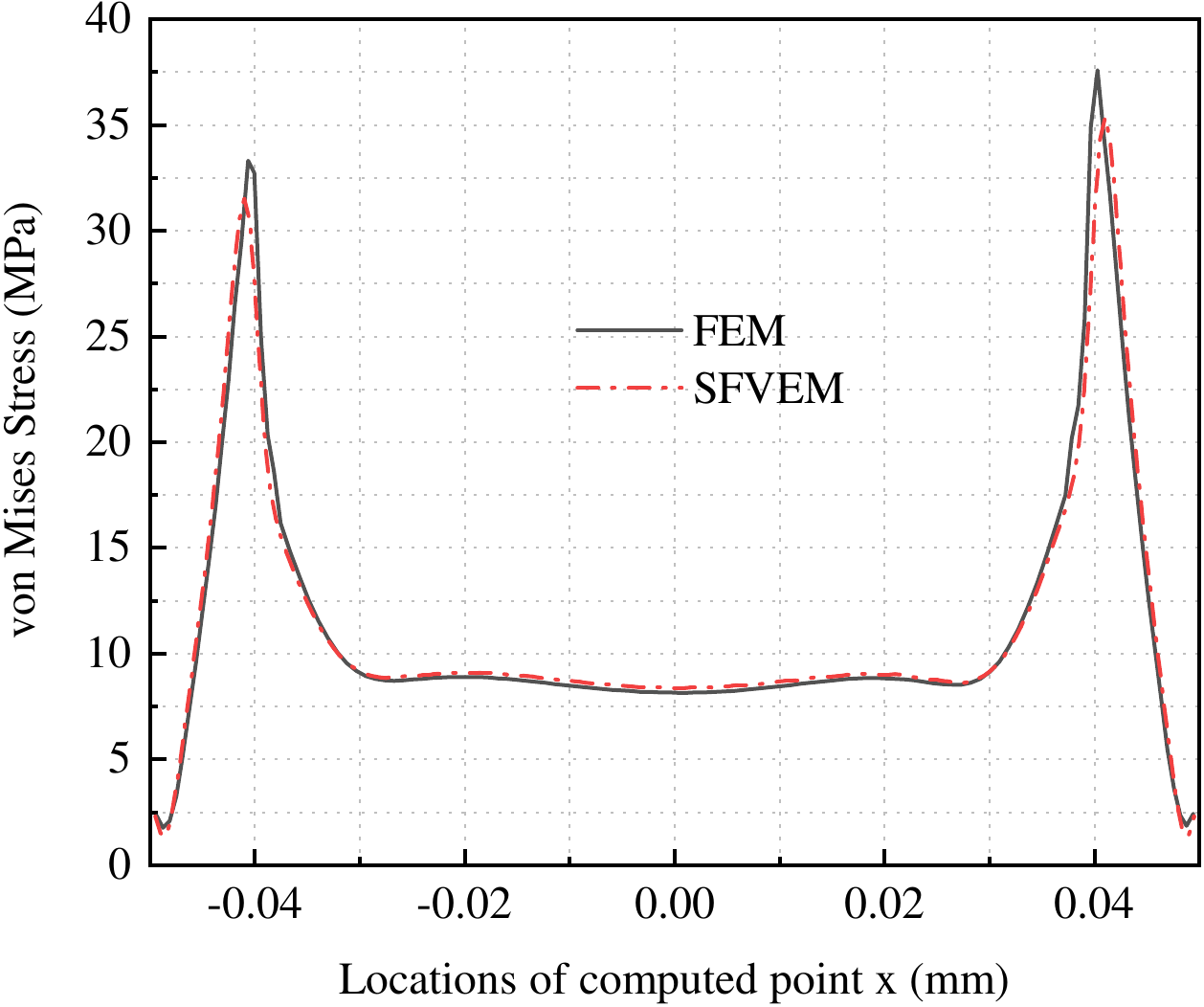}
        \caption{Results along interface $L_1$}
        \label{fig:EX3_von_mises_bgal1}
    \end{subfigure}
     \hfill
    \begin{subfigure}[b]{0.48\textwidth}
        \centering
        \includegraphics[width=0.87\textwidth]{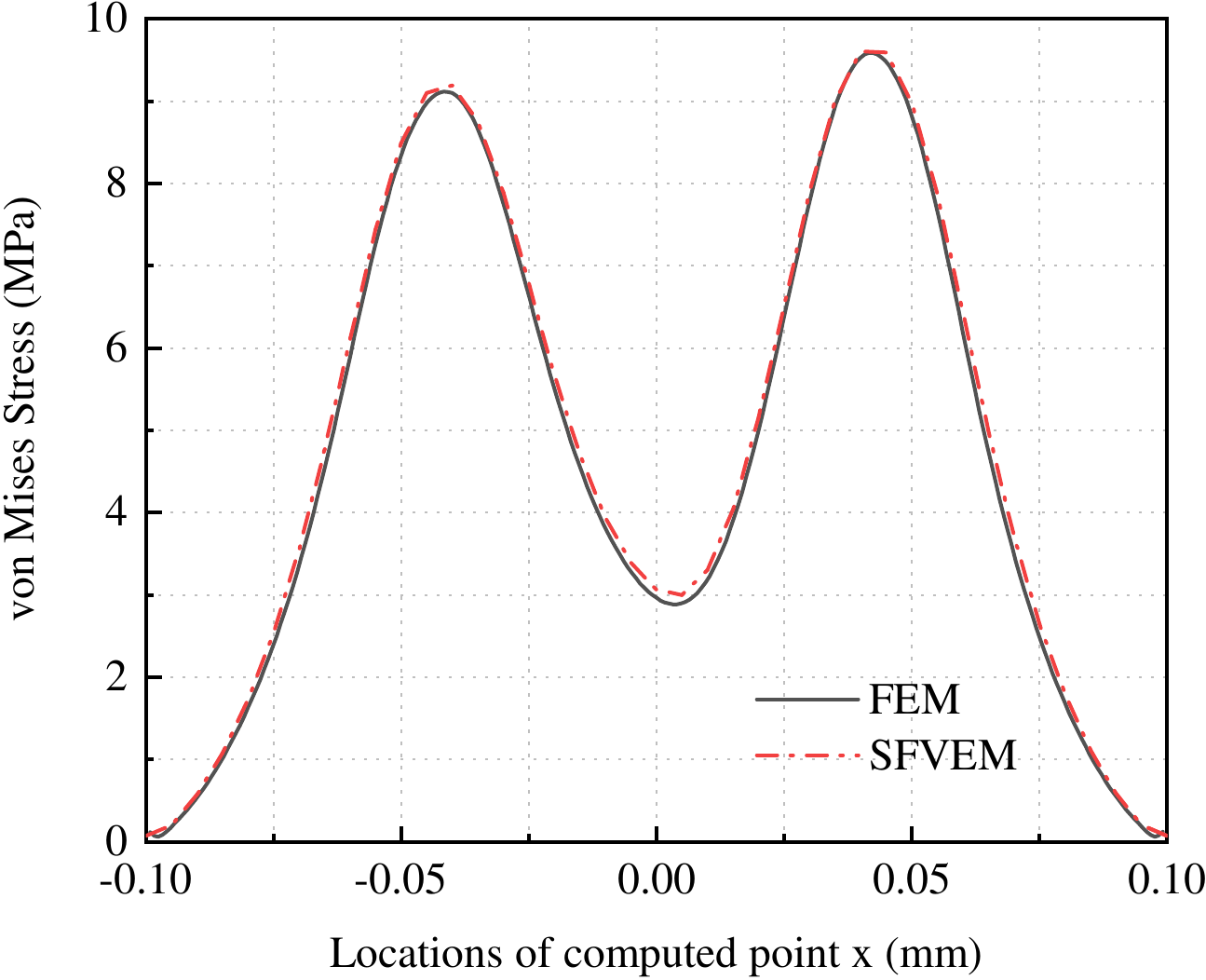}
        \caption{Results along interface $L_2$}
        \label{fig:EX3_von_mises_bgal2}
    \end{subfigure}
     \hfill
    \begin{subfigure}[b]{0.48\textwidth}
        \centering
        \includegraphics[width=0.8\textwidth]{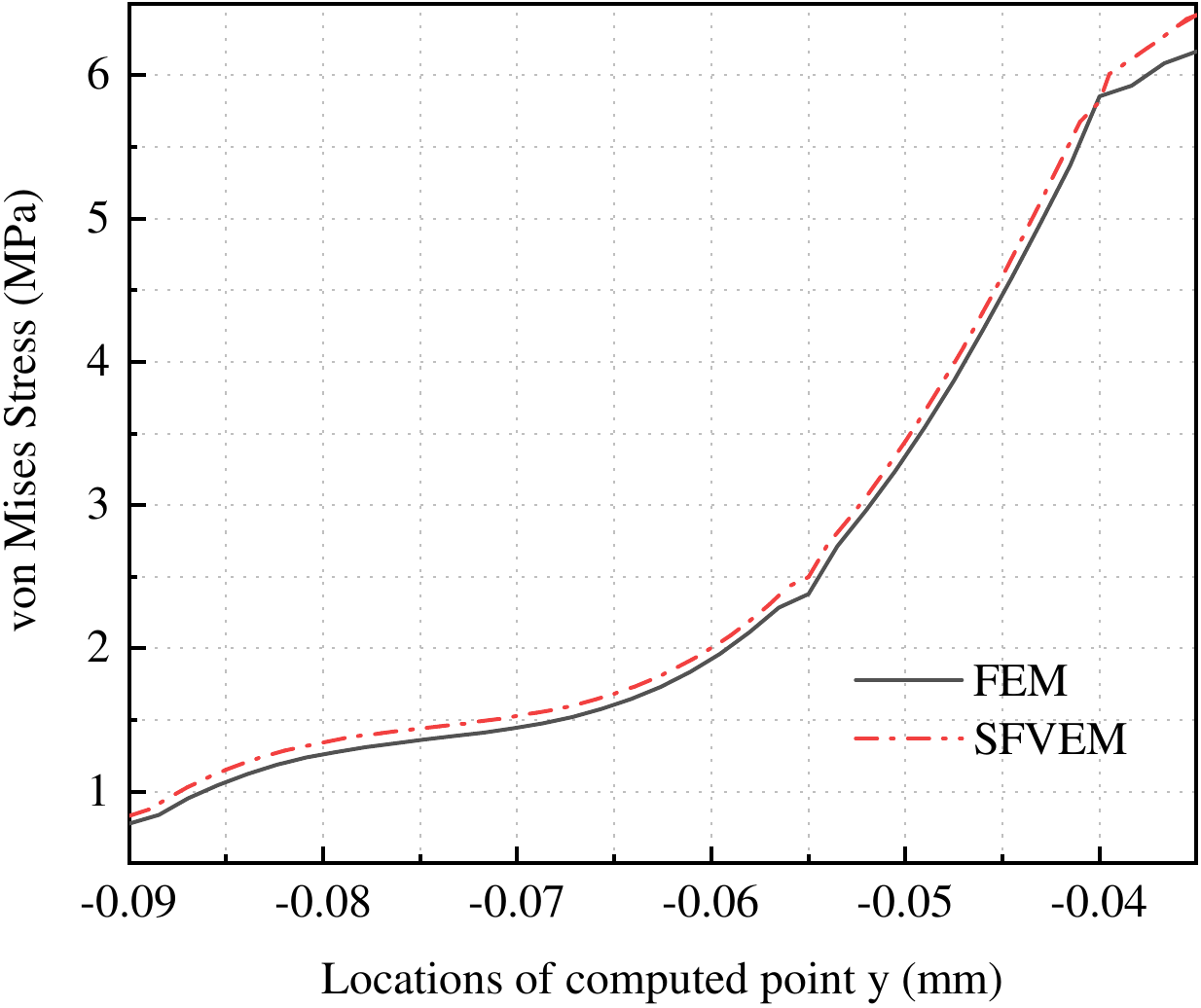}
        \caption{Results along interface $L_3$}
         \label{fig:EX3_von_mises_bgal3}
    \end{subfigure}
     \caption{Stress distribution comparisons along interface $L_1$,$L_2$ and $L_3$.}
   \label{fig:EX3_stress_comparison_bga}
\end{figure}

To validate the proposed SFVEM scheme, we compare the results with a reference FEM solution. 
Figs.~\ref{fig:EX3_von_mises_vem_bga} and~\ref{fig:EX3_von_mises_fem_bga} display the von Mises stress distributions obtained using both methods. The FEM mesh is shown in Fig.~\ref{fig:EX3_bga_mesh_fem}.
The visual comparison demonstrates excellent overall agreement between our SFVEM solution and the reference FEM solution. To further assess the computational accuracy, three lines ($L_1$, $L_2$, $L_3$) are selected in Fig.~\ref{fig:EX3_bga_model}, and the stress distributions along these paths are extracted from the stress distribution contours shown in Fig.~\ref{fig:EX3_stress_comparison_bga1} and compared between FEM and SFVEM. The comparison results are shown in  Figs.~\ref{fig:EX3_von_mises_bgal1}--\ref{fig:EX3_von_mises_bgal3}.


\begin{figure}[htbp]
    \centering
    \includegraphics[width=0.5\textwidth]{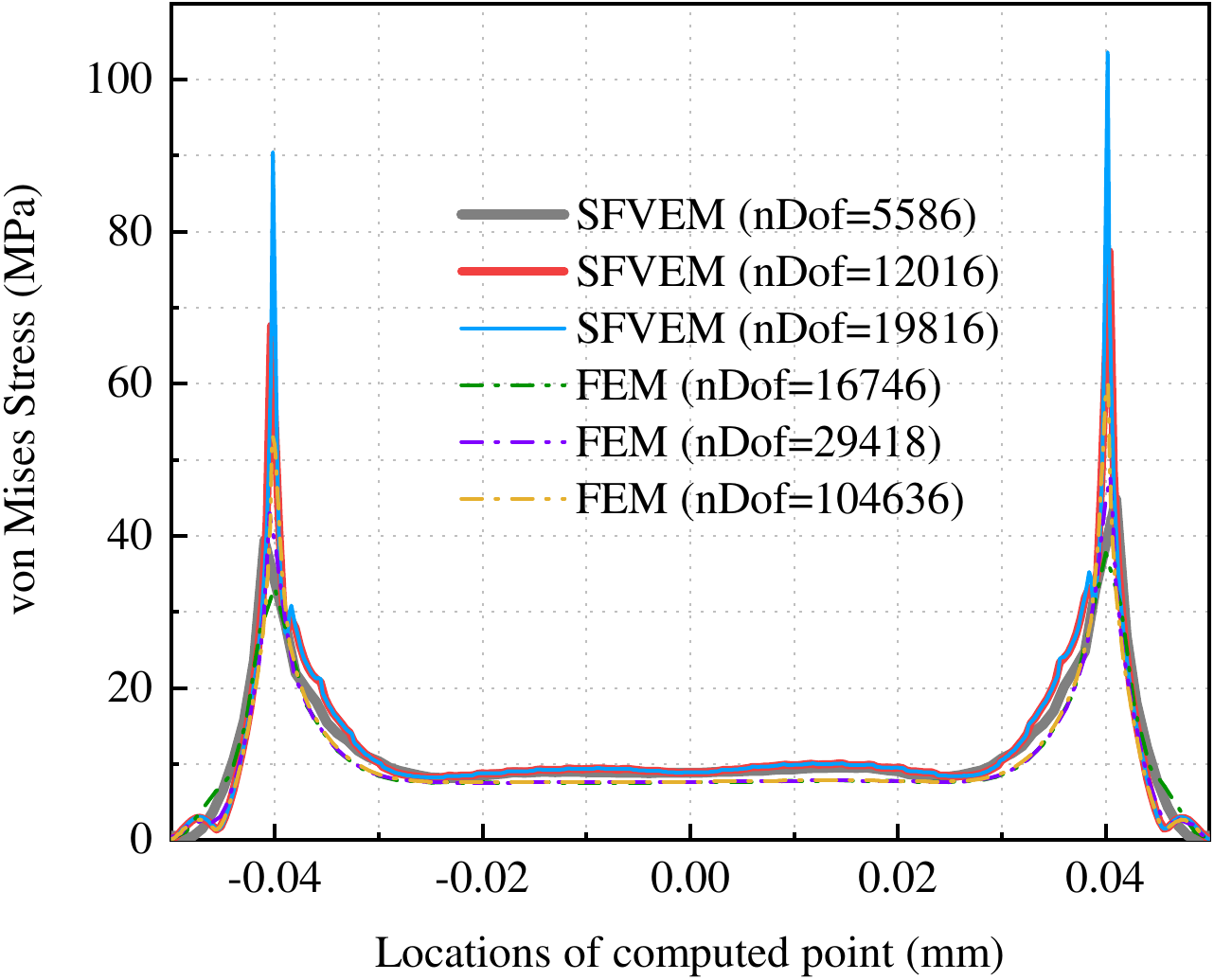}
    \caption{Stress distribution comparison along Al layer upper surface: SFVEM versus locally-refined FEM at multiple discretization levels.}
    \label{fig:EX3_stress_comparison_dof}
\end{figure}

To further illustrate the mesh flexibility and non-matching mesh advantages of our methods, we analyze the stress distribution along the upper surface of the Al layer in the BGA model across different mesh refinement levels. Fig.~\ref{fig:EX3_bga_mesh_fem}  shows the mesh discretization used in the traditional FEM analysis, where smaller element sizes are employed in stress concentration regions while larger elements are used in non-critical areas, with transition elements implemented at coarse-fine mesh interfaces to ensure nodal compatibility.
Fig.~\ref{fig:EX3_stress_comparison_dof} compares results obtained using the SFVEM approach against traditional FEM with local mesh refinement.
The comparative analysis demonstrates remarkable consistency between SFVEM and FEM solutions. 
Notably, even at the lowest mesh density (nDof = 5586), the SFVEM solution maintains excellent agreement with the reference FEM results. 
This confirms the exceptional stability and accuracy of the SFVEM approach even with significantly reduced degrees of freedom, 
highlighting its computational efficiency for complex electronic packaging structures.

\subsection{Structural analysis of plastic ball grid array package}

Interconnect solder joint failure is the most common failure mode in electronic components. 
This failure stems from various factors including temperature cycling, vibration, and mechanical loads. Under real operating conditions, 
printed circuit board assemblies experience not only thermal and vibrational stresses but also mechanical loads such as bending. 
The following analysis examines stress distribution patterns within both the solder ball array and individual solder joints of a PBGA laminated structure subjected to bending loads.

\begin{figure}[htbp]
    \centering
    \includegraphics[width=\textwidth]{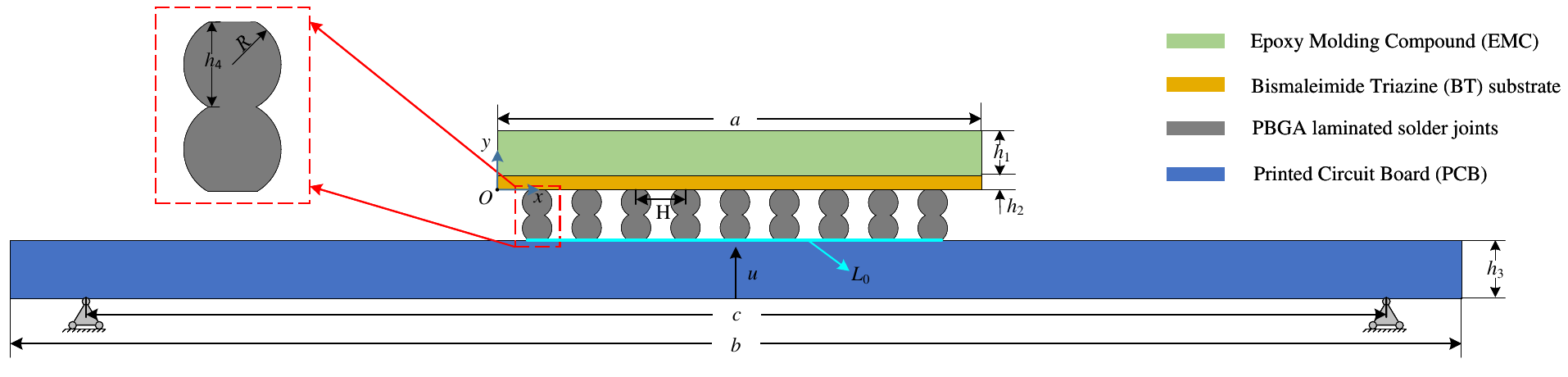}
    \caption{Cross-sectional schematic of PBGA architecture showing layered structure.}
    \label{fig:EX4_pbga_model}
\end{figure}

\begin{table}[htbp]
    \small
    \centering
    \caption{Mechanical properties of PBGA package components.}
    \label{tab:EX4_pbga_properties}
    \begin{tabular}{lcc}
        \toprule
        Material & Young's modulus (MPa) & Poisson's ratio \\
        \midrule
        Molding compound & 155000 & 0.25 \\
        BT substrate & 223000 & 0.30 \\
        PBGA solder joint & 343000 & 0.377 \\
        PCB board & 182000 & 0.25 \\
        \bottomrule
    \end{tabular}
\end{table}

The PBGA laminated solder joint model illustrated in Fig.~\ref{fig:EX4_pbga_model} comprises four key components: epoxy molding compound (EMC), BT substrate, PBGA laminated solder joints, and PCB board. 
The model has the following geometric dimensions: $a = 10\,\text{mm}$, $b = 132\,\text{mm}$, ${h_1} = 0.93\,\text{mm}$, ${h_2} = 0.3\,\text{mm}$, ${h_3} = 1.2\,\text{mm}$, 
with solder balls of diameter $R = 0.61\,\text{mm}$, solder joint height ${h_4} = 0.515\,\text{mm}$, and a solder joint spacing is $H$ = 1 mm. Table~\ref{tab:EX4_pbga_properties} summarizes the material properties for each component. 
The boundary conditions, as shown in Fig.~\ref{fig:EX4_pbga_model}, consist of two fixed displacement constraints at the bottom edges of the PCB board (indicated by the triangular supports, which are symmetrically positioned along the central axis with a spacing of $c = 128\,\text{mm}$.)
while a vertical displacement load of $0.20\,\text{mm}$ is applied upward at the central position of the PCB board's lower surface.

\begin{figure}[htbp]
    \centering
    \includegraphics[width=0.9\textwidth]{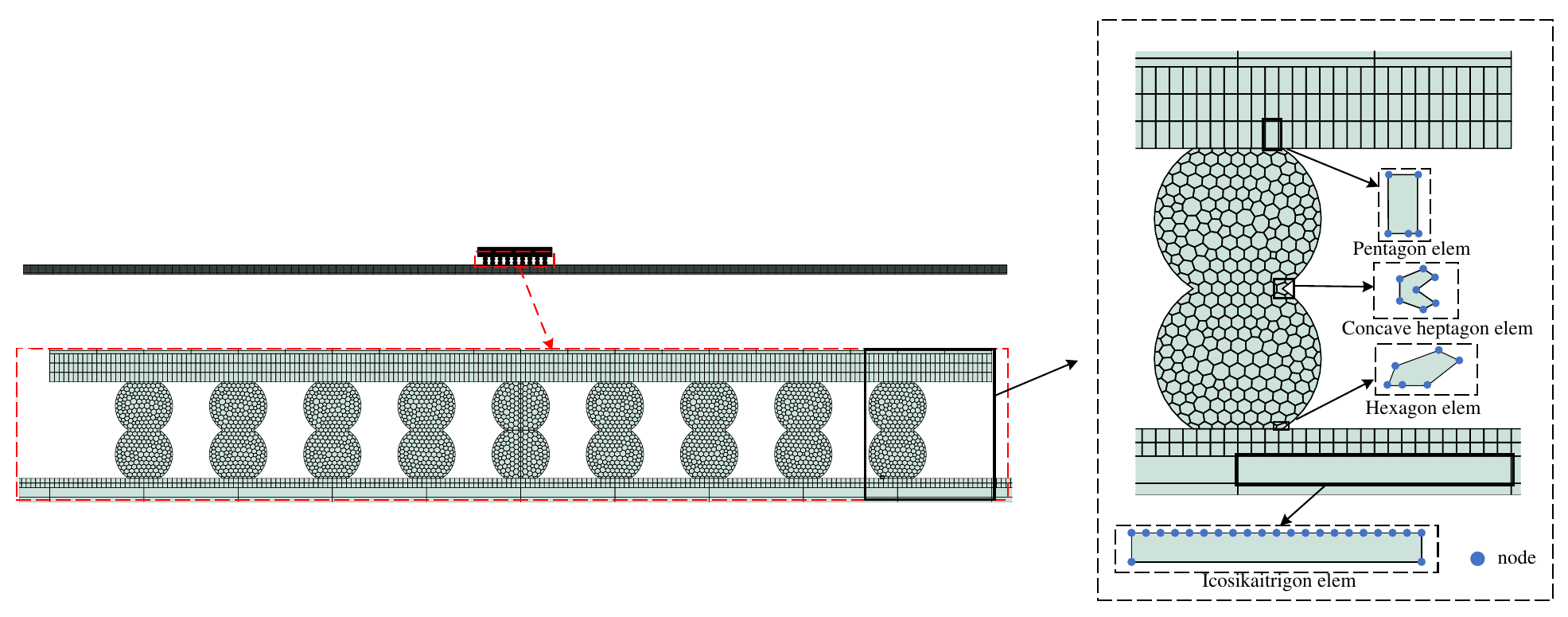}
    \caption{Multi-resolution non-matching mesh for PBGA structure: detailed discretization of solder joints with transitional mesh refinement at material interfaces.}
    \label{fig:EX4_pbga_mesh}
\end{figure}

\begin{figure}[htbp]
    \centering
    \includegraphics[width=0.45\textwidth]{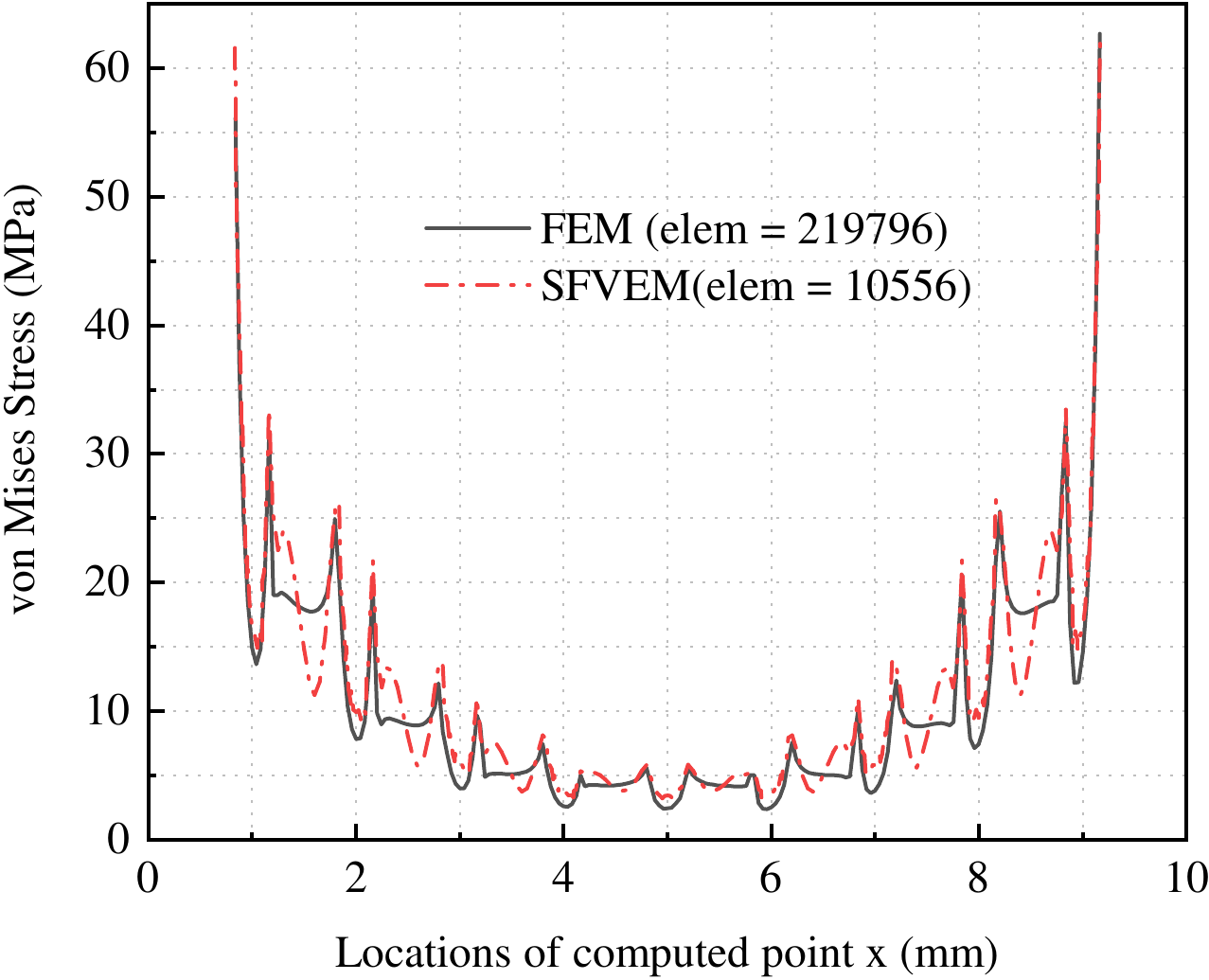}
    \caption{Stress distribution comparison along interface $L_0$: convergence analysis of SFVEM solutions versus reference FEM results at multiple mesh refinement levels.}
    \label{fig:EX4_stress_comparison_L0}
\end{figure}

\begin{figure}[htbp]
    \centering
 \begin{subfigure}[b]{0.9\textwidth}
        \centering
        \includegraphics[width=0.9\textwidth]{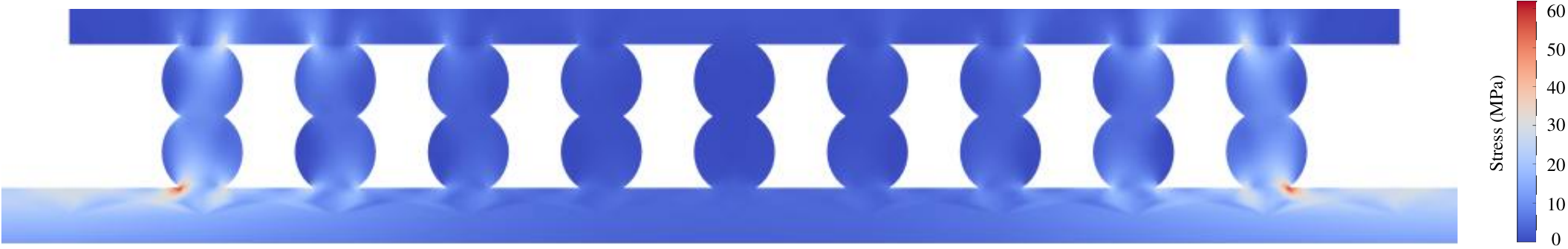}
        \caption{von Mises stress obtained by SFVEM}
        \label{fig:EX4_von_mises_sfvem_PBGA}
    \end{subfigure}
    \hfill
    
    \begin{subfigure}[b]{0.9\textwidth}
        \centering
        \includegraphics[width=0.9\textwidth]{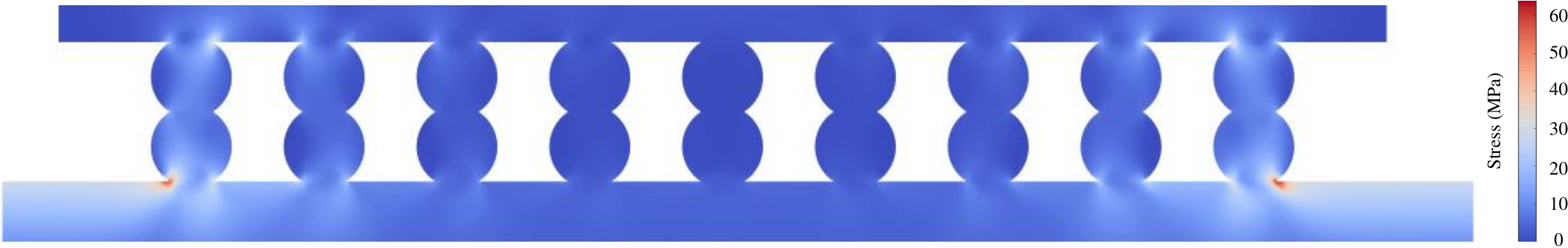}
        \caption{von Mises stress obtained by FEM}
        \label{fig:EX4_von_mises_fem_PBGA}
    \end{subfigure}
    \caption{Comparison of stress distributions around the PBGA laminated solder joints}
    \label{fig:EX1_stress_comparison_PBGA}
\end{figure}
In this model, the PBGA laminated solder joints are relatively small compared to other components. 
FEM analysis requires transition meshes to achieve accurate results, resulting in an excessive number of elements. 
To address this challenge, we implement non-matching meshes as shown in Fig.~\ref{fig:EX4_pbga_mesh}.

We perform comparative analysis using both the proposed SFVEM and FEM, extracting the stress distribution at the interface $L_0$ labeled in Fig.~\ref{fig:EX4_pbga_model}.
Fig.~\ref{fig:EX4_stress_comparison_L0} illustrates these stress distributions, demonstrating excellent consistency between the SFVEM and FEM results using refined meshes. Figs.~\ref{fig:EX4_von_mises_sfvem_PBGA} and~\ref{fig:EX4_von_mises_fem_PBGA} display the von Mises stress distributions around the PBGA laminated solder joints obtained using both methods.
The analysis confirms that maximum stress in the solder ball array occurs at the outermost solder joints, which aligns with established failure patterns in PBGA packages under bending loads.

\begin{figure}[htbp]
    \centering
    \includegraphics[width=0.15\textwidth]{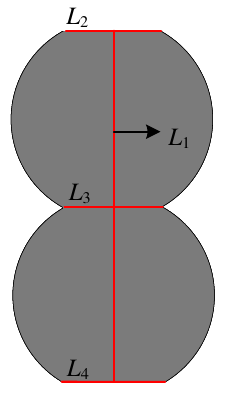}
    \caption{Central solder joint detail: critical interfaces and measurement locations for stress analysis.}
    \label{fig:EX4_pbga_central_solder}
\end{figure}

\begin{figure}[htbp]
    \centering
    \begin{subfigure}[b]{0.48\textwidth}
        \centering
        \includegraphics[width=0.9\textwidth]{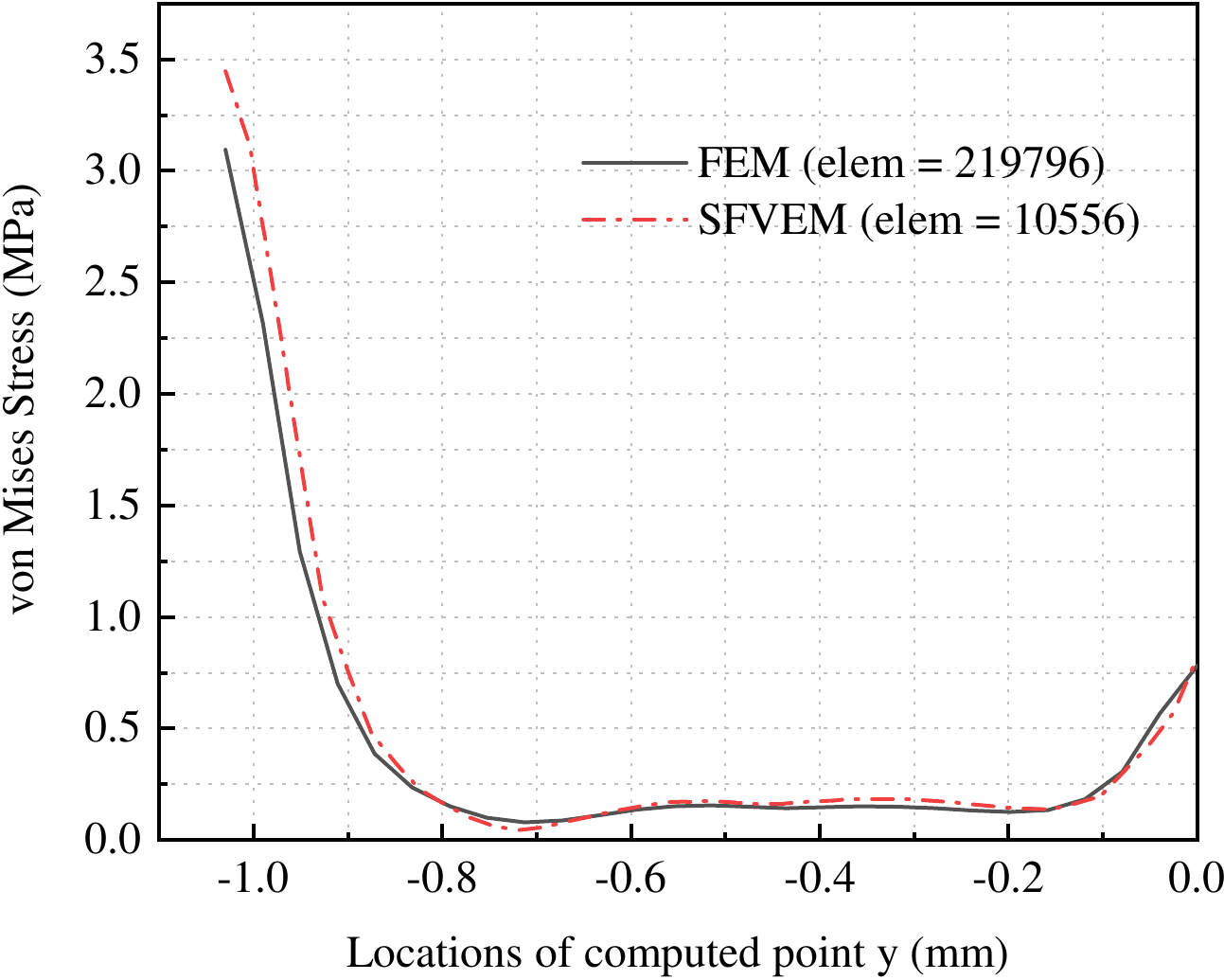}
        \caption{Results along interface $L_1$}
        \label{fig:EX4_central_stress_b}
    \end{subfigure}
    \hfill
    \begin{subfigure}[b]{0.48\textwidth}
        \centering
        \includegraphics[width=0.9\textwidth]{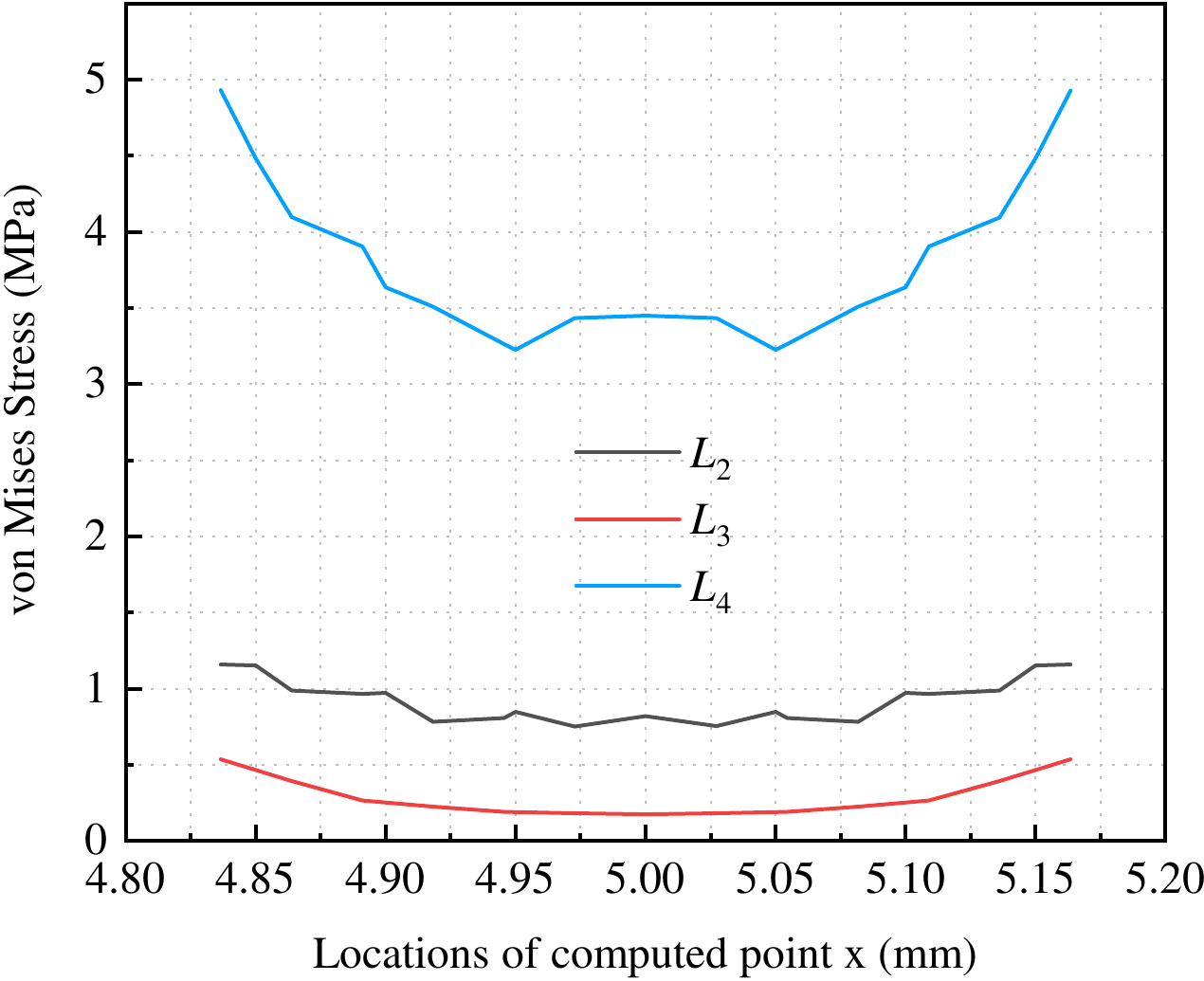}
        \caption{Results along interfaces $L_2$, $L_3$, and $L_4$}
        \label{fig:EX4_central_stress_a}
    \end{subfigure}
    \caption{Stress distribution analysis within central solder joint: (a) comparative SFVEM-FEM results along interface $L_1$; (b) stress variations across multiple interfaces.}
    \label{fig:EX4_central_joint_stress}
\end{figure}

We further analyze stress distributions at multiple interfaces ($L_1$, $L_2$, $L_3$, $L_4$) within the central solder joint, as shown in Fig.~\ref{fig:EX4_pbga_central_solder}. 
Fig.~\ref{fig:EX4_central_stress_a} presents the stress distribution profiles along interfaces $L_2$, $L_3$, and $L_4$. 
Due to both structural and loading symmetry, the stress distributions exhibit symmetric patterns about $x = 5$, with maximum stress occurring at the bottom interface $L_4$. 
Fig.~\ref{fig:EX4_central_stress_b} compares SFVEM and FEM results along interface $L_1$, 
demonstrating excellent agreement between the methods and confirming that the SFVEM accurately captures the same stress variation patterns as the reference FEM solutions. The minor differences observed, particularly at lower $y$ values, are primarily attributed to the different mesh densities employed: SFVEM uses 10,556 elements compared to FEM's 219,796 elements. This computational efficiency advantage of SFVEM comes with the expected trade-off in local resolution, though convergence studies confirm that SFVEM results approach the FEM reference solution with mesh refinement.

\subsection{Elastic analysis of representative volume element (RVE) sub-model}

\begin{figure}[htbp]
    \centering
    \begin{subfigure}[b]{0.48\textwidth}
        \centering
        \includegraphics[width=0.5\textwidth]{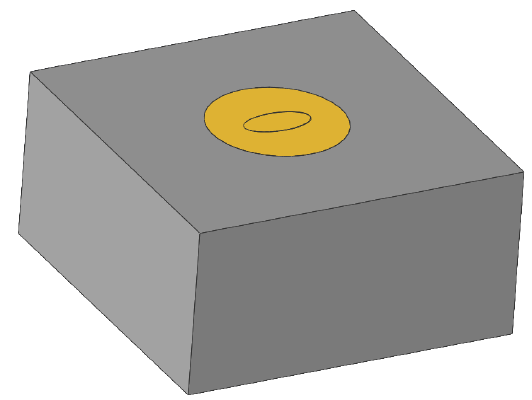}
        \caption{Cross-sectional RVE geometry}
        \label{fig:EX5_rve_model}
    \end{subfigure}
    \hfill
    \begin{subfigure}[b]{0.48\textwidth}
        \centering
        \includegraphics[width=0.5\textwidth]{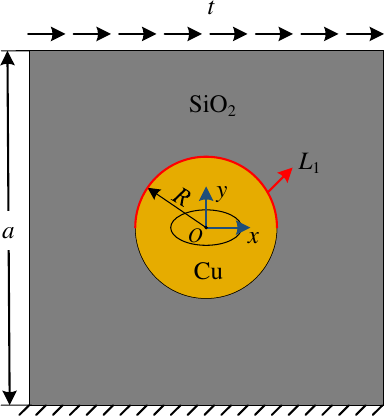}
        \caption{Planar configuration (top view)}
        \label{fig:EX5_rve_top_view}
    \end{subfigure}
    \caption{TGV RVE: (a) structural cross-section and (b) top surface configuration.}
    \label{fig:EX5_rve_diagrams}
\end{figure}

Fig.~\ref{fig:EX5_rve_model} illustrates the Representative Volume Element (RVE) sub-model, with particular focus on the two-dimensional top surface. 
The model has the following geometric parameters: $a = 1\,\text{mm}$, copper via diameter $R = 0.2\,\text{mm}$, with elliptical features having a major axis of $0.2\,\text{mm}$ and minor axis of $0.1\,\text{mm}$. 
The material properties in our simulation include: SiO$_2$ substrate with Young's modulus $E_{\text{SiO}_2} = 75\,\text{GPa}$ and Poisson's ratio $\nu_{\text{SiO}_2} = 0.17$; 
Copper interconnects with Young's modulus $E_{\text{Cu}} = 150\,\text{GPa}$ and Poisson's ratio $\nu_{\text{Cu}} = 0.3$. 
For boundary conditions, we apply fixed constraints at the bottom surface and a uniform stress of $2\,\text{MPa}$ on the top surface. It should be noted that while this geometry represents an RVE configuration, the primary objective of this example is to demonstrate the mesh flexibility and robustness of the SFVEM method rather than to perform periodic analysis. Therefore, periodic boundary conditions are not implemented; instead, the left and right boundaries are treated as free (traction-free) boundary conditions, as shown in Fig.~\ref{fig:EX5_rve_top_view}.
\begin{figure}[htbp]
    \centering
    \begin{subfigure}[b]{0.48\textwidth}
        \centering
        \includegraphics[width=0.85\textwidth]{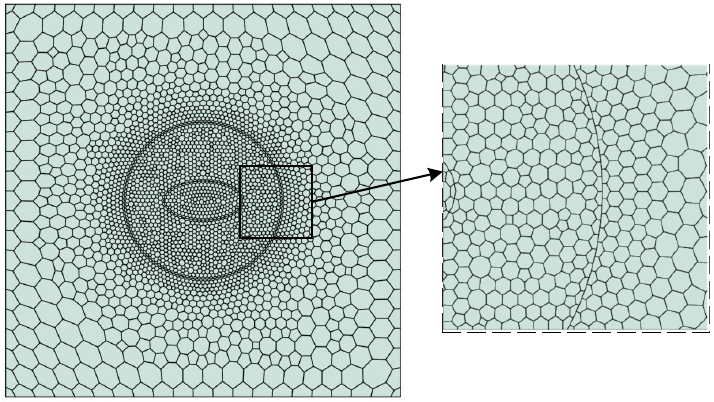}
        \caption{$\theta = 0^{\circ}$ orientation}
        \label{fig:EX5_mesh_0deg}
    \end{subfigure}
    \hfill
    \begin{subfigure}[b]{0.48\textwidth}
        \centering
        \includegraphics[width=0.85\textwidth]{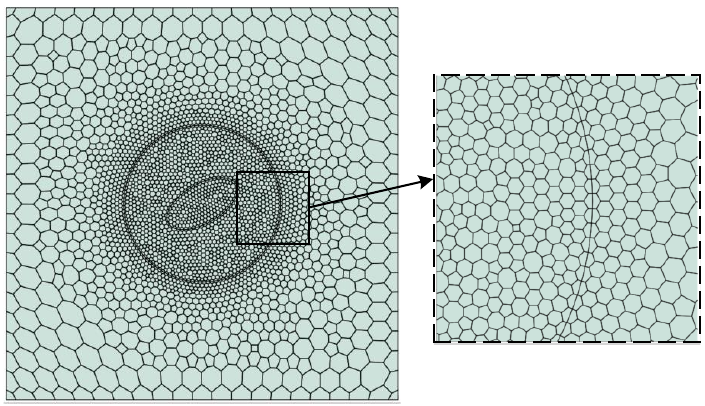}
        \caption{$\theta = 30^{\circ}$ orientation}
        \label{fig:EX5_mesh_30deg}
    \end{subfigure}
    \\ 
    \begin{subfigure}[b]{0.48\textwidth}
        \centering
        \includegraphics[width=0.85\textwidth]{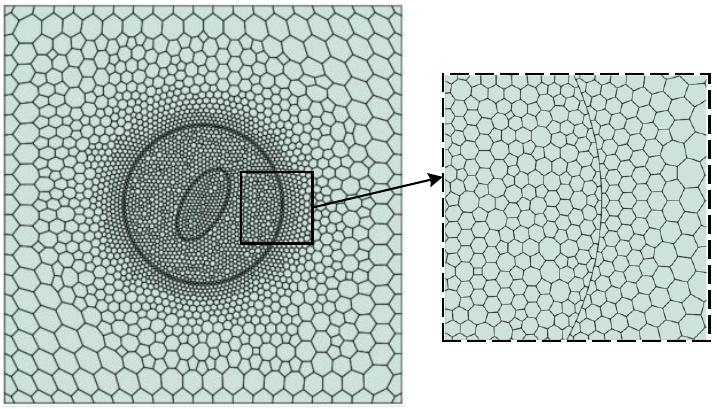}
        \caption{$\theta = 60^{\circ}$ orientation}
        \label{fig:EX5_mesh_60deg}
    \end{subfigure}
    \hfill
    \begin{subfigure}[b]{0.48\textwidth}
        \centering
        \includegraphics[width=0.85\textwidth]{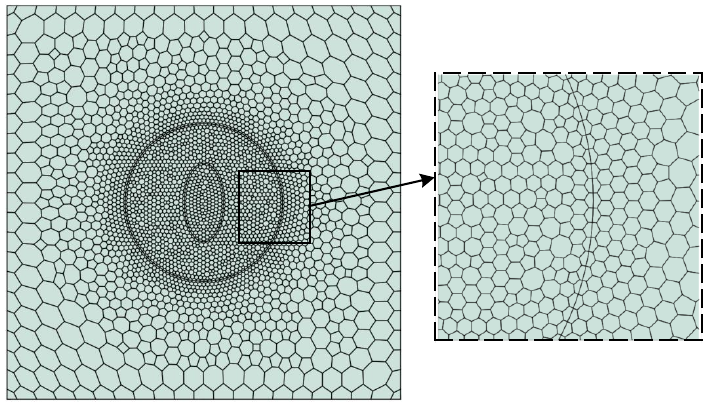}
        \caption{$\theta = 90^{\circ}$ orientation}
        \label{fig:EX5_mesh_90deg}
    \end{subfigure}
    \caption{Copper component mesh configurations at various rotational orientations (measured from horizontal axis).}
    \label{fig:EX5_mesh_rotations}
\end{figure}

In the analysis of this model, the mesh discretization is shown in Figs.~\ref{fig:EX5_mesh_0deg}--\ref{fig:EX5_mesh_90deg}. 
Figs. \ref{fig:EX5_mesh_30deg}, \ref{fig:EX5_mesh_60deg}, and \ref{fig:EX5_mesh_90deg} display the copper component mesh rotated by 30°, 60°, and 90° respectively from the original orientation in Fig.~\ref{fig:EX5_mesh_30deg}. 
The elliptical feature serves as a visual indicator of the mesh rotation angle within the copper domain.
As shown in the figure, the size, shape, and number of elements within the model remain fixed and unchanged. 
The only modification is the rotation angle of the mesh within the Cu region. Non-matching meshes are implemented at the Cu-SiO$_2$ interface, with a total element count of 3,041.

\begin{figure}[htbp]
    \centering
    \includegraphics[width=\textwidth]{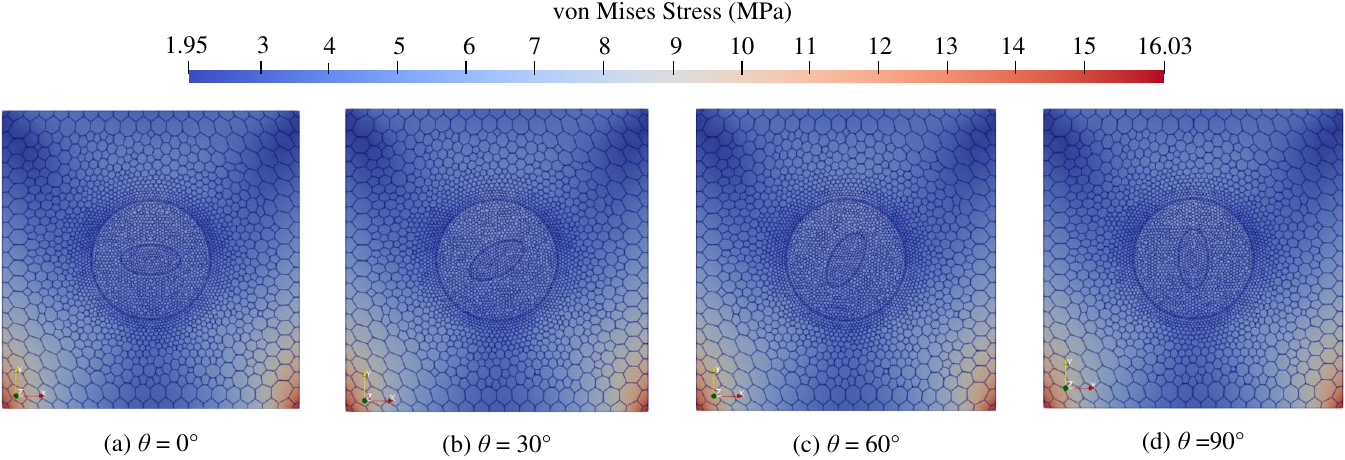}
    \caption{Comparative stress distribution across copper component with varying mesh orientations ($\theta$).}
    \label{fig:EX5_stress_comparison_rotations}
\end{figure}

\begin{figure}[htbp]
    \centering
    \includegraphics[width=0.5\textwidth]{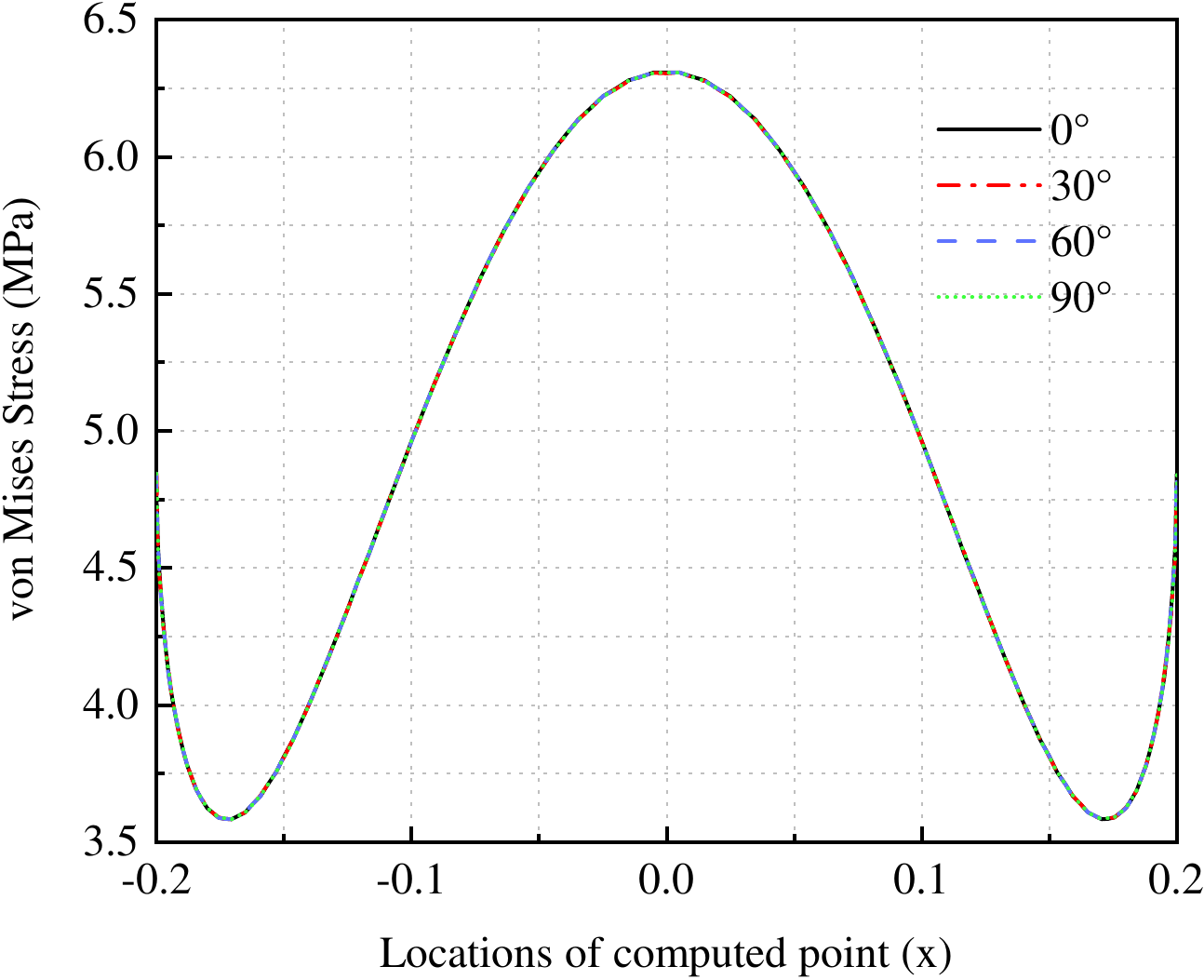}
    \caption{The stress distribution along the arc $L_1$.}
    \label{fig:EX5_RVE stress}
\end{figure}

Fig.~\ref{fig:EX5_stress_comparison_rotations} presents the stress contours obtained from all four mesh configurations. 
Fig.~\ref{fig:EX5_RVE stress} shows the stress distribution along the upper surface of the inner circle $($ the arc $L_1$ shown in Fig.~\ref{fig:EX5_rve_top_view} $)$ of the SiO$_2$ component for different mesh rotation angles. The results clearly demonstrate that both stress distribution patterns and magnitude values remain consistent across all mesh orientations. 
This consistency further validates the exceptional mesh flexibility of our SFVEM algorithm and highlights the significant advantage of non-matching mesh implementation, 
where localized mesh modifications can be performed without affecting the overall mesh structure or solution accuracy.

\section{Conclusions}
\label{conclusions}
This work presents a Virtual Element Method for analyzing electronic packaging structures with the following key contributions: 
(1) We developed a novel stabilization-free virtual element method for thermomechanical problems.
The stabilization-free virtual element method does not require any additional stabilization terms, so the calculation results are not affected by parameters, which improves the reliability of the results.
(2) Through comprehensive numerical examples, we rigorously verified the proposed method's accuracy and effectiveness.
Our method demonstrates excellent correlation with analytical solutions in both heat conduction and thermoelastic analyses. 
(3) Our method shows significant advantages when analyzing complex electronic packaging structures, including TSV, BGA, and PBGA models. 
VEM effectively manages stress concentrations and interface conditions, producing reliable results for both thermal and mechanical responses. 
(4) We demonstrated VEM's inherent mesh flexibility across various configurations. 
Our results confirm that solution accuracy remains consistent regardless of mesh rotation or element shape, validating the method's robustness for practical applications.

This methodology provides an efficient framework for analyzing electronic packaging structure reliability, particularly for problems involving multiple geometric scales and complex geometries. The complete source code and implementation details for the algorithms presented in this work are available on our GitHub repository (https://github.com/yanpeng-gong/VEM-electronic-packaging) and the VEMhub website (www.vemhub.com). Future research will focus on extending this SFVEM approach to three-dimensional thermomechanical analysis of electronic packaging structures, building upon existing 3D SFVEM frameworks, incorporating more sophisticated material behaviors, and investigating fully coupled solutions for cases where strong bidirectional coupling effects are present.

\section*{Acknowledgements}
\label{sec:Acknowledgements}
The research was supported by the National Natural Science Foundation of China (No. 12002009). 

\section*{Appendix A: Analytical Solution Derivation for Thick-Walled Cylinder}

For the one-dimensional, steady-state temperature field equation without internal heat generation
$$\frac{\mathrm{d}}{\mathrm{d} r}\left(r \frac{\mathrm{d} T}{\mathrm{~d} r}\right)=0$$

With boundary conditions $T=T_a$ at $r=a$ and $T=T_b$ at $r=b$, where $a$ and $b$ are the inner and outer radii respectively. Integrating both sides yields
$$r \frac{\mathrm{d} T}{\mathrm{~d} r}=c_1 \Rightarrow T=c_1 \ln r+c_2$$

Applying boundary conditions
$$\left\{\begin{array}{l}T_a=c_1 \ln a+c_2 \\ T_b=c_1 \ln b+c_2\end{array}\right.$$

Solving for the unknown parameters
$$c_1=\frac{T_b-T_a}{\ln (b / a)}, \quad c_2=T_a-\frac{T_b-T_a}{\ln (b / a)} \ln a$$

The analytical temperature solution becomes:
$$T=T_a+\frac{T_b-T_a}{\ln (b / a)} \ln (r / a)$$

For the mechanical problem in polar coordinates, the strain-displacement relationships are
$$\varepsilon_r=\frac{\mathrm{d} u_r}{\mathrm{~d} r}, \quad \varepsilon_\theta=\frac{u_r}{r}$$

Under plane stress conditions, the stress-strain relationships are
$$\left\{\begin{array}{l}\sigma_r=\frac{E}{1-v^2}\left(\varepsilon_r+v \varepsilon_\theta-\alpha(1+v) T(r)\right) \\ \sigma_\theta=\frac{E}{1-v^2}\left(\varepsilon_\theta+v \varepsilon_r-\alpha(1+v) T(r)\right)\end{array}\right.$$

The equilibrium equation in polar coordinates (neglecting body forces) is
$$\frac{\partial \sigma_r}{\partial r}+\frac{1}{r}\left(\sigma_r-\sigma_\theta\right)=0$$

Substituting the stress-strain relationships into the equilibrium equation yields
$$\frac{\mathrm{d}^2 u_r}{\mathrm{~d} r^2}+\frac{1}{r} \frac{\mathrm{~d} u_r}{\mathrm{~d} r}-\frac{1}{r^2} u_r=\alpha(1+v) \frac{T_b-T_a}{r \ln (b / a)}$$

This is an Euler equation with both homogeneous and particular solutions. The homogeneous solution has the form $u^1(r)=B_1 r+B_2 r^{-1}$, and the particular solution is $u^2(r)=D r \ln r$ with
$$D=\alpha(1+v) \frac{T_b-T_a}{2 \ln (b / a)}$$

The complete displacement solution is
$$u(r)=B_1 r+B_2 r^{-1}+D r \ln r$$

The stress components become
$$\left\{\begin{array}{l}\sigma_r=\frac{E}{1-v^2}\left(B_1(1+v)+B_2(v-1) r^{-2}+D(\ln r(v+1)+1)-\alpha(1+v) T(r)\right) \\ \sigma_\theta=\frac{E}{1-v^2}\left(B_1(1+v)+B_2(1-v) r^{-2}+D(\ln r(v+1)+v)-\alpha(1+v) T(r)\right)\end{array}\right.$$

Applying stress-free boundary conditions $\sigma_r(r=a)=\sigma_r(r=b)=0$ and solving the resulting system of equations determines the constants $B_1$ and $B_2$ as given in Eqs.\eqref{eq:coefficient_B1} - \eqref{eq:coefficient_B2} .

\bibliography{VEM_ref} 

\end{document}